\documentclass{article}
\usepackage{amssymb}
\usepackage{graphicx}
\usepackage{amsmath}

\newtheorem{theorem}{Theorem}[section]

\newtheorem{conjecture}[theorem]{Conjecture}
\newtheorem{corollary}[theorem]{Corollary}

\newtheorem{definition}[theorem]{Definition}

\newtheorem{lemma}[theorem]{Lemma}

\newtheorem{proposition}[theorem]{Proposition}
\newtheorem{remark}[theorem]{Remark}

\newcommand{\SL}{\mathrm{SL}}
\newcommand{\GL}{\mathrm{GL}}
\newcommand{\EL}{\mathrm{EL}}

\newcommand{\Hom}{\mathrm{Hom}}
\newcommand{\Mat}{\mathrm{Mat}}
\newcommand{\mb}[1]{\mbox{\rm {#1}}}
\newcommand{\mr}[1]{\mathcal{#1}}
\newcommand{\Rd}{\left.R^*\right.^2}
\newcommand{\la}{\langle}
\newcommand{\ra}{\rangle}
\newcommand{\C}{\mathbb{C}}

\newcommand{\R}{\mathbb{R}}
\newcommand{\T}{\mathbb{T}}
\newcommand{\Z}{\mathbb{Z}}
\newcommand{\N}{\mathbb{N}}
\newcommand{\F}{\mathbb{F}}
\newcommand{\Id}{\mathrm{Id}}
\newcommand{\KC}{\mr{K}}
\newcommand{\TC}{\tau}

\newcommand{\TauC}{$\tau$-constant}

\newcommand{\KaC}{Kazhdan constant}

\begin{document}
\title{Universal lattices and unbounded rank expanders }
\author{M. Kassabov } \date{}
\author{Martin Kassabov}

\maketitle
{
\renewcommand{\thefootnote}{}
}
\begin{abstract}
We study the representations of non-commutative universal lattices
and use them to compute lower bounds for the
\TauC\ for the commutative universal lattices
$G_{d,k}= \SL_d(\Z[x_1,\dots,x_k])$ with respect to
several generating sets.

As an application of the above result we show that the
Cayley graphs of the finite groups $\SL_{3k}(\F_p)$ can be made
expanders using suitable choice of the generators.
This provides the first examples of expander families of groups
of Lie type where the rank is not bounded and gives a
natural (and explicit) counter examples to two
conjectures of Alex Lubotzky and Benjamin Weiss.
\end{abstract}

{
\renewcommand{\thetheorem}{\arabic{theorem}}
\section*{Introduction}

The groups $\SL_d(\mr{O})$, where $\mr{O}$ is a ring of integers is
in a number field $K$, have many common --- they all
have Kazhdan property \emph{T}, have a finite congruence kernel and super rigidity.
In~\cite{YSh},
Y.~Shalom conjectured that many of these properties are inherited form
the group $\SL_d(\Z[x])$. He called the groups
$\SL_d(\Z[x_1,\dots,x_k])$ \emph{universal lattices}.
Almost nothing is
known about the representation theory of these groups.
It is conjectured that the universal lattices have property \emph{T},
it was known that these groups have property
$\tau$, see~\cite{KNtau}.

In this paper we study the noncommutative analogs of these groups, i.e., the subgroups
of $\GL_d(\Z\la x_1,\dots, x_k\ra)$ generated by all elementary matrices
with coefficients in the free associative ring on $k$ generators
$\Z\la x_1,\dots, x_k\ra$. Using that these groups can be mapped onto
many universal lattices we obtain new bounds for the \TauC\ for the groups
$\SL_d(\Z[x_1,\dots,x_k])$ with respect to several generating sets.

\medskip

Let $G$ be a topological group and consider the space
$\widetilde{G}$ of all equivalence classes of unitary
representations of $G$ on some Hilbert space $\mr{H}$. This space
has a naturally defined topology, called the Fell topology, as
explained in \cite{LubZuk}  \S1.1 or~\cite{expanderbook}
Chapter 3 for example. Let $1_G$ denote the
trivial $1$-dimensional representation of $G$ and let
$\widetilde{G}_0$ be the set of representations in $\widetilde{G}$
which do not contain $1_G$ as a subrepresentation
(i.e., do not have invariant vectors).

\begin{definition}
A group is $G$ is said to have Kazhdan property
\emph{T} if $1_G$ is isolated from
$\widetilde{G}_0$ in the Fell topology of $\widetilde{G}$.
\end{definition}

A discrete group with property \emph{T} is finitely generated: see
\cite{kazhdan}. Since in this paper we shall be concerned mostly
with discrete groups we give the following equivalent
reformulation of property \emph{T}:

\emph{Equivalent definition:}
Let $G$ be a discrete group generated by a finite set $S$.
Then $G$ has the Kazhdan property \emph{T} if there is
$\epsilon  >0$
such that for every nontrivial irreducible unitary
representation $\rho: \ G \rightarrow U(\cal{H})$ on a
Hilbert space $\mr{H}$ and every  vector $v \not = 0$ there
is some $s \in S$ such that $||\rho(s)v-v||> \epsilon ||v||$.
The largest $\epsilon$ with this property is called the \emph{\KaC}
for $S$ and is denoted by $\mr{K}(G;S)$.

\medskip

The property \emph{T} depends only on the group $G$
and does not depend on the choice of the
generating set $S$, however the Kazhdan constant depends also on the generating set.

Property \emph{T} implies certain group theoretic conditions on $G$
(finite generation, FP, FAB etc) and can be used for construction
of \emph{expanders} from the finite images of $G$.
For this last application the following weaker property $\tau$
is sufficient:

Let $\widetilde{G}^f$ and $\widetilde{G}^f_0$ denote the finite representations
of $\widetilde{G}$, resp. $\widetilde{G}_0$ (i.e., the representations which factor
through a finite index subgroup).

\begin{definition}
A group is $G$ is said to have property $\tau$ if $1_G$
is isolated from $\widetilde{G}^f_0$
in the induced Fell topology of $\widetilde{G}^f$.
Equivalently: the group $G$ with the
profinite topology and its profinite completion $\widehat{G}$
has property \emph{T}.
\end{definition}

Again, for a discrete finitely generated%
\footnote{A discrete group with property $\tau$ may not be finitely generated. However, all the groups
we are going to work with are f.g., and there is no need to define \TauC\ for not finitely generated groups.}
group $G$ there
is an equivalent definition of these properties:

\emph{Equivalent definition:}
Let $G$ be an discrete group generated by a set $S$.
Then $G$ has the property $\tau$  if there is
$\epsilon  >0$
such that for every nontrivial finite
irreducible  unitary
representation $\rho: \ G \rightarrow U(\cal{H})$ on a
Hilbert space $\mr{H}$ and every  vector $v \not = 0$ there
is some $s \in S$ such that $||\rho(s)v-v||> \epsilon ||v||$.
The largest $\epsilon$ with this property is called the \emph{\TauC}
and is denoted by $\TC(G;S)$.

Property $\tau$ is not interesting for groups which do not have
many finite quotients. All the groups we are going to work with are
residually finite and have enough finite factors.

\medskip

As mentioned above, properties $T$ and $\tau$
can be used for the
construction of families of \emph{expanding graphs}.
We need more terminology. \medskip

A finite graph $\Gamma$ is called an $\epsilon$-expander for some
$\epsilon \in (0,1)$ if for any subset $A \subseteq \Gamma$ of size
at most $|\Gamma|/2$ we have $|\partial (A)| >\epsilon |A|$
(where $\partial(A)$ is the set of vertices of $\Gamma \backslash A$ of edge
distance 1 to $A$). The largest such $\epsilon$ is called the
expanding constant of $\Gamma$. Constructing families of $\epsilon$-expanders
with large expanding constant $\epsilon$ and bounded valency
is an important practical problem in computer science. For an excellent
introduction to the subject we refer the reader to the
book~\cite{expanderbook} by Lubotzky.
The following equivalent definition of
property $\tau$ illustrates one of the few
known approaches to explicit construction of expanders:
\begin{proposition}%
[\cite{expanderbook}, Theorem 4.3.2]
\label{CGexpanders}
Let $G$ be a finitely generated discrete group. Then $G$ has
property $\tau$ iff  for any set of generators $S$ there exists an
$\epsilon=\epsilon(S) >0$ such that
the Cayley graphs $\mr{C}(G_i,S_i)$ of the finite images of
$G_i$ of $G$ (with respect to the images $S_i$ of $S$)
form a family of $\epsilon$-expanders. The largest $\epsilon_0(S)$
with this property is related to the \TauC\ $\TC(G,S)$ defined above.
In particular we have
$\epsilon_0(S) \geq \TC(G,S)^2/4$.
\end{proposition}

\medskip

Our approach to property $\tau$ is inspired by a paper by
Shalom \cite{YSh} which relates property \emph{T} to
\emph{bounded generation property} of high rank Chevalley groups
over rings of integers. In this paper we will work only with
the group $\EL_d(R)$ for $d\geq 3$ which correspond to the Dynkin
diagram $A_{d-1}$. Similar arguments also work for high rank the Chevalley
groups which arise from others Dynkin diagrams.

Let $R$ be an associative ring with unit, and for $i \not = j \in \{1,2,...,d\}$ let
$E_{i,j}$ denote the set of elementary $d \times d$ matrices
$\{\Id+r\cdot e_{i,j}\  | \ r \in R \}$. Set
$E=E(R)=\bigcup_{i\not =j}E_{i,j}$ and
let $\EL(d;R)=\EL_d(R)$ be the subgroup of the multiplicative group of the ring of
$d\times d$ matrices over $R$ generated by $E(R)$.

\begin{definition}
The group $G=\EL_d(R)$ is said to have
\emph{bounded elementary generation property}
if there is a
number $N=BE_d(R)$ such that every element of $G$
can be written as a product of at most $N$
elements from the set $E$.
\end{definition}

Examples of infinite rings $R$ satisfying the above definition are rings of
integers $\mr O$ in number fields $K$ (for $d\geq 3$), see \cite{CK}.
In this classical case this property is
known as bounded generation because each group $E_{i,j}\simeq (\mr O,+)$
is a product of finitely many cyclic groups.

\medskip

The following theorem is a generalization of
a result of Y.~Shalom (see Theorem 3.4 in~\cite{YSh})
to non-commutative rings:

\begin{theorem}
\label{t1}
Suppose that $d\geq 3$ and $R$ is a finitely generated associative ring
such that $\EL_d(R)$ has bounded elementary generation property.
Then $\EL_d(R)$ has property \emph{T} (as a discrete group).

Moreover there is an explicit lower bound for the \KaC\ $\KC(G,S)$
with respect to a specific generating set $S$.
\end{theorem}

\bigskip

A very interesting conjecture is whether
$
\EL_d(\Z[x])=\SL_d(\Z[x])
$
and other universal lattices,
has bounded elementary generation property.
This together with Theorem~\ref{t1} implies that
the universal lattices has property \emph{T}.
In~\cite{KNtau} it was shown that the groups
$$
G_{d,k}:=
\EL_d(\Z[x_1,\dots,x_k])=\SL_d(\Z[x_1,\dots,x_k])
$$
have property $\tau$.
We believe that the methods in
this paper can be used to show that the universal lattice in the
non-commutative case have property $\tau$. In order to show that
we need several results about the structure of finite non-commutative ring $R$
and the $K$-groups $K_1(R)$ and $K_2(R)$.

Theorem~\ref{t1} has several applications:

\subsection*{Tau Constant for $\SL_d(\Z[x_1,\dots, x_k])$.}

In~\cite{KNtau}, the commutative universal lattices are shown to
have property $\tau$ and
the \TauC\ for $G_{d,k}=\SL_d(\Z[x_1,\dots, x_k])$ is estimated.
Theorem~5 from~\cite{KNtau} gives a lower bound for the \TauC\ with respect to
the generating set $S_{d,k}$, which is
asymptotically $O(d^{-2}22^{-k})$.

However, it is possible to improve this estimate to $O(d^{-1/2}(1+k/d)^{-3/2})$.
We state this result as Theorem~\ref{taucom}, its proof
is an extension of the one in~\cite{KNtau} and combines the ideas from~\cite{K},
namely relative property \emph{T} of certain groups and the
use of \emph{generalized elementary matrices},
see Theorem~\ref{relTpq} and Lemma~\ref{congimagest}.

\begin{theorem}
\label{taucom}
The commutative universal lattices
$$
G_{d,k}= \EL_d(\Z[x_1,\dots,x_k]) = \SL_d(\Z[x_1,\dots,x_k])
$$
have property $\tau$. The \TauC\ with respect to the generating set $S_{d,k}$
consisting of all elementary matrices with $\pm 1$ of the diagonal and the ones
with $\pm x_i$ next to the main diagonal, satisfies
$$
\TC(G_{d,k},S_{d,k}) \geq
\frac{1}{800\sqrt{d}\left(1 +(k/d)^{3/2} \right)}.
$$
\end{theorem}

It is interesting if this bound is asymptotically exact.
Using a natural representation of $\SL_d(\Z[x_1,\dots,x_k])$
in the Hilbert space $l^2(\F_2^d)$ it can be seen that
$$
\TC(G_{d,k},S_{d,k}) \leq \sqrt{\frac{2}{d}},
$$
which shows that if $k$ is fixed then the bound for \TauC\
in Theorem~\ref{taucom} is asymptotically exact.
We are not aware of any nontrivial upper bound for $\TC(G_{d,k},S_{d,k})$
which depends on the number of generators $k$.

If one modifies the generating set $S_{d,k}$ then the \TauC\ can be improved,
in fact Theorem~\ref{tconsforlargerings} gives that for many $d$ there exist a
finite generating set $S_{d,k}'$ such that
$\TC(G_{d,k},S_{d,k}')$ depends only on $k/d$.

\subsection*{Unbounded rank expanders.}

Theorem~\ref{t1} has another interesting application.
In~\cite{lubwiess} Lubotzky and Weiss made the following conjectures:

\begin{conjecture}[\cite{lubwiess}, Problem~5.1]
\label{conjexp}
Let $\Gamma_i$ be a family of finite groups and let $\Sigma_i$ and $\Sigma_i'$ be two
generating sets of the group $\Gamma_i$ of uniformly bounded size.
suppose that the Cayley graphs $\mr{C}(\Gamma_i,\Sigma_i)$ form an expander family. Is it true that
that $\mr{C}(\Gamma_i,\Sigma_i')$ is also an expander family?
\end{conjecture}

\begin{conjecture}[\cite{lubwiess}, Conjecture~5.4]
\label{conjambT}
Let $K$ be a compact group. If $\Gamma_1$ and $\Gamma_2$ are two finitely generated
dense subgroups of $K$ with $\Gamma_1$ amenable and $\Gamma_2$ having property $T$
then the group $K$ is finite.
\end{conjecture}

In~\cite{alonlub} Alon, Lubotzky and Wigderson showed that
Conjecture~\ref{conjexp} is false using zig-zag products. Their construction
is based on a probabilistic arguments and does not give explicit generating
sets which make the Cayley graphs expanders.
Using Theorem~\ref{t1}, we are able to construct more natural
counter examples of this conjecture, moreover in our counter example all groups and
generating sets are explicit and the expanding constants can be estimated.

Consider the family of groups $\Gamma_{n,p} = \SL_n(\F_p), n>3$.
Let $\Sigma_n$ be the generating set of $\Gamma_{n,p}$ consisting of
$$
A_n = \sum_{i=1}^{n-1} e_{i,i+1} + (-1)^{n-1}e_{n,1} \quad \mbox{and} \quad B_n = \Id_n+e_{1,2}
$$
and their inverses. In~\cite{lubwiess}, it is shown that the \KaC s of these groups are bounded form above by
$$
\KC(\Gamma_{n,p},\Sigma_n) \leq \sqrt{2/n},
$$
which implies that any infinite subset of the Cayley graphs $\mr{C}(\Gamma_{n,p},\Sigma_n)$
for fixed $p$ and different $n$, are not expanders.%
\footnote{It is interesting to note that although these Cayley are not expanders, they
have a logarithmic diameter, see~\cite{KR}}
In view of Conjecture~\ref{conjexp} this implies that
any infinite family of the Cayley graphs of $\Gamma_{n,p}$ are not expanders with respect to any
generating set of uniformly bounded size.

\medskip

Using non-commutative universal lattices $G^n_{g,k}$ it is possible to
find generating sets $\Sigma'_{3l}$ of $\Gamma_{3l,p}$ such that the
Caylay graphs are expanders:

Let us observe that the matrix algebras $\Mat_l(\F_p)$ can be
generated as rings by $2$ elements, i.e., all these algebras
 are quotients
of the ring $\Z\la x_1,x_2\ra$. It is well known that
$$
\Mat_d(\Mat_l(R)) \simeq \Mat_{dl}(R).
$$
These two observations allow us to consider the groups
$\SL_{3l}(\F_p)$ as quotients of $\EL_3(\Z\la x_1,x_2\ra)$.
Theorem~\ref{t1} gives us that the group
$\EL_3(\Z\la x_1,x_2\ra)$ almost has property $\tau$ with respect
to the generating set $S_{3,2}$ consisting of $28$ elements.
If we project this generating set to a generating set $\Sigma'_{3l}$
of $\SL_{3l}(\F_p)$ we obtain that \KaC \ is uniformly bounded from
bellow and therefore the Cayley graphs $\mr{C}(\SL_{3l}(\F_p), \Sigma_{3l}')$
form an expander family.%
\footnote{This argument can be generalized to show that for any $n\geq 3$ and any prime $p$,
there exists a generating set $\Sigma''_{n,p}$ of $\SL_{n}(\F_p)$, such that
the Cayley graphs $\mr{C}(\SL_{n}(\F_p), \Sigma_{n,p}'')$ for all $n$ and $p$
form an expander family.}
In fact we have shown that
\begin{theorem}
\label{UbRExpanders}
For all primes $p$ and $l\geq 1$,
there exists a generating set $\Sigma'_{3l}$ of the group $\SL_{3l}(\F_p)$
with $28$ elements
such that the \KaC s
$$
\KC(\SL_{3l}(\F_p), \Sigma'_{3l}) > {1}/{400} \quad \mbox{for all } l,p,
$$
which implies that the Cayley graphs
$\mr{C}(\SL_{3l}(\F_p), S'_{3l})$ form a family of expanders
with expanding constant at least $1.5 \times 10^{-6}$.
\end{theorem}

Theorem~\ref{UbRExpanders} gives a positive answer to the question~\cite[Problem 5.27]{LubZuk},
whether the Cayley graphs of $\SL_n(\F_p)$ for fixed $p$ can be made expanders.

\medskip

Similar construction almost allows us to disprove Conjecture~\ref{conjambT} --- it is known
(see~\cite{lubwiess}) that in the product $\prod_n \SL_{n}(\F_p)$
there is a finitely generated dense amenable subgroup. We can show that
there is also a dense subgroup which has property $\tau$.%
\footnote{It is possible that
this group also have property \emph{T}, but we are unable to prove it.}
In order to simplify the argument we will construct a counter example inside
the product $\prod_l \SL_{3l}(\F_p)$:
\begin{theorem}
\label{amtau}
Let $p>2$ be a prime and let $\tilde K^p$ be the infinite compact group
 given by the product:
$$
\tilde K^p:= \prod_{l\geq 3} \SL_{3l}(\F_p)
$$
with the the product topology. There
exist finitely generated dense subgroups $\Gamma_1$ and $\Gamma_2$ of $K$ such that
$\Gamma_1$ is amenable and $\Gamma_2$ has property $\tau$.
\end{theorem}

\bigskip

This idea can be applied to a slightly more general situation.
There is a surjection $R_{s+2}^n \to \Mat_l(R^n_{sl})$, which gives a
surjection $G_{3,s+2}^n \to G_{3l,sl}^n$. Using this surjection we can
prove that:
\begin{theorem}
\label{tconsforlargerings}
There exists a set $S_{3l,sl}'$ with $28 + 8 s$ elements generating the
group $G_{3l,sl}=\SL_{3l}(\Z[x_1,\dots,x_{sl}])$ such that the \TauC\
$
\TC(G_{3l,sl}, S_{3l,sl}')
$
depends only on $s$.
\end{theorem}

\bigskip

A few words about the structure of the rest of the paper:
Section~\ref{outline} contains the proof of Theorems~\ref{t1} and~\ref{taucom} modulo Lemmas~\ref{abe}
and~\ref{congimage}. In section~\ref{relative} we prove some technical results about
relative property \emph{T} and relative \KaC s. In section~\ref{uni} we prove
uniform bounded elementary generation of the group $\SL_d(\bar R)$ for some finite $\bar R$.
The last section~\ref{applications} is devoted to applications of Theorem~\ref{t1} --
a construction of expander family of Lie groups of unbounded rank and a counter examples of the
conjectures of A.~Lubotzky and B.~Wiess.


\subsection*{Acknowledgements}
I am grateful to  R.~K.~Dennis, A.~Lubotzky, N.~Nikolov, Y.~Shalom and L.~Vaserstain for
introducing me to the subject and for useful discussions. I additionally wish to thank
Nikolay Nikolov for providing the the proof of Lemma~\ref{congimage} and others
results in section~\ref{uni}.

\subsection*{Notations}

Let $R^n_k:= \Z\la x_1,...,x_k\ra$ be the free associative ring on $k$ (non-commutative)
generators, and let $R_k:= R_k^n/[R_k^n,R_k^n] = \Z[x_1,...,x_k]$ be its abelinization, which
is isomorphic to the polynomial ring on $k$ variables with coefficients in $\Z$.

Let $R$ be an associative ring with $1$, consider the
ring $\Mat_d(R)$ of $d\times d$ matrices of $R$.
With $\GL_d(R)$ we will denote the multiplicative group $\Mat_d(R)$, i.e., the group
of invertible $d\times d$ matrices with coefficients in the ring $R$. We will assume that $d \geq 3$.

For $i \not = j \in \{1,2,...,d\}$ let
$E_{i,j}$ denote the set of elementary $d \times d$ matrices of the from
$\{\Id+r\cdot e_{i,j}\  | \ r \in R \}$, it is clear that these matrices are invertible.
Also set
$E=E(R)=\bigcup_{i\not =j}E_{i,j} \subset \GL_d(R)$
and let $\EL(d;R)=\EL_d(R)$
be the subgroup
in $\GL_d(R)$ generated by $E(R)$.
In the case of finitely generated commutative ring $R$ this is a subgroup%
of  $\SL_d(R)$.
For the rest of this paper we will denote
$G^n_{d,k}:=\EL(d;R_k^n)$ and $G_{d,k}:=\EL(d;R_k)$.
Unless we explicitly need to specify the number of generators of the ring
$R_k^n$ or the size of matrices $G_{d,k}^n$ and $G_{d,k}$
we will denote $R_k^n$ by $R^n$, $R_k$ by $R$, $G_{d,k}^n$ by $G^n$ and $G_{d,k}$ by $G$.

The groups $G_{d,k}^n$ and $G_{d,k}$ are generated%
\footnote{The sets $S_{d,k}$ generate the groups $G_{d,k}$ and $G^n_{d,k}$ provided that $d\geq 3$.
If $d=2$ these groups are not finitely generated.} by
the set $S_{d,k}= F_1\cup F_2$, where
$F_1$ is the set of $2(d^2-d)$ elementary matrices with $\pm 1$ off the diagonal
and $F_2$ is the set $4k(d-1)$ elementary matrices $\Id \pm x_le_{ij}$ with
$|i-j|=1$ and $1\leq l \leq k$.

For any $r\in R$ and $1\leq i\not=j\leq d$,
with $e_{ij}(r)$ we will denote the matrix $\Id + r\cdot e_{ij} \in \EL_d(R) \subset \GL_d(R)$.
If $r$ is a unit of the ring $R$, $e_{ii}(r)$ will denote the matrix $\Id + (r-1)\cdot e_{ii}
\in \GL_d(R)$.

We will call the elements in $E(R)$ elementary matrices.
Any matrix in $\EL_{p+q}(R)$ of the form
$$
\left(
\begin{array}{cc}
\Id_p & * \\ 0 & \Id_q
\end{array}
\right)
\quad \mb{or} \quad
\left(
\begin{array}{cc}
\Id_q & 0 \\ * & \Id_p
\end{array}
\right)
$$
where $0$ and $*$ are blocks of sizes $p\times q$ and $q\times p$ respectively is called a
\emph{generalized elementary matrix} (abbreviated to GEM).

We will also assume that all rings are associative with $1$.
}

\section{Proof of Theorems~\ref{t1} and~\ref{taucom}}
\label{outline}

\textbf{Proof of Theorem~\ref{t1}:}
Let $\rho : \EL_d(R) \to U(\mr{H})$ be a unitary representation of the group
$\EL_d(R)$ for some $k$-generated associative ring $R$, which have the
bounded elementary generation property. Suppose that $v\in \mr{H}$ is an
$\epsilon$-almost invariant unit vector for the set $S$, where
$\epsilon \leq  K_{d,k}$, here $K_{d,k}$ is a constant which will be determined later.

The next lemma is a corollary to Theorems~\ref{relT} and~\ref{relTpq},
proved in Section \ref{relative}.
\begin{lemma}
\label{abe}
There exist constants
$$
M(k)<3\sqrt{2} (\sqrt{k}+3)\quad \mbox{ and }\quad M_d(k) < 6\sqrt{2}(\sqrt{k}+3)+\sqrt{3d}
$$
with the following properties:
Let $\rho$ be a unitary representation of the group
$G_{d,k}^n=\EL(d;R_k^n)$ in a Hilbert space $\mr{H}$.
Let $v\in \mr{H}$ be a unit vector such that
$||\rho(g)v-v|| <\epsilon$ for $g\in S_{d,k}$. Then:

a) We have that $||\rho(g)v - v || \leq 2M(k)\epsilon$
for every elementary matrix $g$.

b) We have that
$||\rho(g)v - v || \leq 2M_d(k)\epsilon$
for every generalized elementary matrix $g$.
\end{lemma}

\bigskip

To finish the proof of Theorem~\ref{t1} we use an argument due to
Y.~Shalom from~\cite{YSh}.
Let $g$ be any element in the group $\EL_d(R)$. Using the bounded elementary generation property of $R$, we
can write
$$
g=\prod_{s=1}^{BE_d(R)} u_s,
$$
where $u_s$ are elementary matrices.
Finally, we can use Lemma~\ref{abe} to obtain:
$$
||\rho(g)v-v|| \leq \sum_{s=1}^{BE_d(R)} ||\rho(u_s)v-v|| <
\sum_{s=1}^{BE_d(R)} 2 M(k)\epsilon  = 2BE_d(R) M(k)\epsilon.
$$
If we chose $K_{d,k} = \left(\sqrt{2}BE_d(R)M(k)\right)^{-1}$ and use that $\epsilon \leq  K_{d,k}$ we obtain
$$
||\rho(g)v-v|| < 2 BE_d(R) M(k) K_{d,k} = \sqrt{2}.
$$
This show that every element in the group $\EL_d(R)$ moves the vector $v$ by less then $\sqrt{2}$.
Let $V$ denote the the $G^n$-orbit of the vector $v$ in $\mr{H}$. The center of mass $c_V$
of the set $V$ is invariant under the action of $G$ and is not zero because the orbit $V$
lies entirely in the half space $\{v'\in \mr{H} \mid \Re (v,v') > 0\}$.
Therefore $\mr{H}$ contains a nontrivial $G$-invariant vector, which completes the
proof that the group $\EL_d(R)$ has the property \emph{T} and also provides a
lower bound for the \KaC. $\square$

\bigskip

\noindent
\textbf{Proof of Theorem~\ref{taucom}:}
We can not apply Theorem~\ref{t1} directly because
it is not known whether the group
$G_{d,k} = \SL_d(\Z[x_1,\dots,x_k])$ has bounded elementary generation
property. We want only to estimate the \TauC\ for this group, not the \KaC, and
it is sufficient to show that all finite images of $G_{d,k}$ have uniform
bounded elementary generation property%
.

Let $\rho : G_{d,k} \to U(\mr{H})$ be a unitary representation of the
universal lattice $G_{d,k}$ which factors through a
finite index subgroup. Suppose that $v\in \mr{H}$ is an $\epsilon$-almost
invariant unit vector for the set $S_{d,k}$, where
$\epsilon \leq  T_{d,k}$. Here $T_{d,k}$ is a constant depending on $d$ and $k$
which will be fixed later.

Let $H = \ker \rho < G$, this is a normal subgroup of finite index in $G$.
Define subgroups $H_{i,j}$ of $H$ and subsets $U_{i,j}$ of the ring $R$ for $i \not =j$ as follows
$$
H_{i,j}:= H \cap E_{i,j}; \quad \quad U_{i,j}:=\{r\in R | \,\,\, 1 + r\cdot e_{i,j} \in H\}.
$$
Using elements in the Weyl group we can see that the subgroups $E_{i,j}$ are pairwise conjugate in $\SL_d(R)$.
This gives us that the sets $U_{i,j}$ do not depend on the indices $i,j$.
Using commutation with a suitable elementary matrix we can see that $U_{i,j}$ is also closed under
multiplication by elements in $R$
, i.e., $U=U_{i,j}$ is an ideal
of the ring $R$. This ideal is of finite index in $R$ because $H$ is a subgroup of finite index in $G$.

Let $\EL(d;U)$ be the normal subgroup in $G$ generated by $H_{i,j}$:
$$
\EL(d;U):= \langle H_{i,j} \rangle^G < H
$$

Let $SK_1(R,U;d)=G(U)/\EL(d;U)$ where $G(U)$ is the principal congruence subgroup of
$G$ modulo $U$, i.e., the matrices in $G$ congruent to $\Id_d$ modulo $U$.
We have the following diagram
$$
\begin{array}{c@{}c@{}c@{}c@{}c@{}c@{}c@{}c@{}c}
1 & \rightarrow & SK_1(R,U;d) & \rightarrow & G/\EL(d;U) & \rightarrow & G(R/U)= \SL_d(\bar{R}) & \rightarrow & 1 \\
  &             &             &             & \downarrow&             &                        &             & \\
  &             &             &             &   G/H     &             &                        &             & \\
  &             &             &             & \downarrow&             &                        &             & \\
  &             &             &             &    1      &             &                        &             & \\
\end{array}
$$
where the row and the column are exact.

\medskip

The next result is known in the case of commutative rings
(cf. \cite{DV}, Remark 10 in the commutative case). Its proof
in case of matrices over finite commutative ring is slightly more complicated:
\begin{lemma}
\label{congimage}
Let $\bar R = \Mat_t(R')$ be a matrix ring over a finite commutative ring $R'$.
The group $\EL_d(\bar{R}) = \SL_{dt}(R')$ has uniform bounded elementary generation property, i.e.,
number $BE_d(\bar R)$ of elementary matrices needed to express any element in
$\EL_d(\bar R)$ is depends only on $d$ and is independent of the ring $\bar R$.

In fact every matrix in $\EL_d(\bar{R})$ can be written as a product of $38$ generalized
elementary matrices.
\end{lemma}

\medskip

It remains to deal with the group $SK_1(R,U;d)$. It is always a finite abelian group and
it measures the departure of $\EL(d;U)$ from being
a congruence subgroup.
The following lemma is proved in~\cite{KNtau}:
\begin{lemma}[\cite{KNtau}, Lemma~10]
\label{kernel}
Each element of the groups $K_2(R/U)$ is a \linebreak[4]
product of at most
$k+2$ Steinberg symbols $\{a,b\}$,
for $a,b$ invertible elements in $R/U$. The
same is true for its image  $SK_1(R,U;d)$.
\end{lemma}

For a $u, u' \in R$ such that $uu' = 1 \mod U$ and a pair of indices $i,j$ consider the elements
$$
w_{i,j}(u,u'):=e_{i,j}(u)e_{j,i}(-u')e_{i,j}(u) \textrm{ and } h_{i,j}(u,u'):=w_{i,j}(u,u')w_{i,j}(-1,-1),
$$
in $\SL_d(R)/\EL(d;U)$. It is clear that the elements $w_{i,j}(u,u')$ and $h_{i,j}(u,u')$
depend only on the image of $u$ in $R/U$ and we can denote them $w_{i,j}(\bar u)$ and $h_{i,j}(\bar u)$.
For any units $\bar u, \bar v \in R/U$ we define the Steinberg symbol
$$
\{\bar u,\bar v\}_{i,j}=h_{i,j}(\bar u\bar v) h_{i,j}(\bar u)^{-1}h_{i,j}(\bar v)^{-1} \in \SL_d(R)/\EL(d;U).
$$
It can be shown that the element $\{\bar u,\bar v\}_{i,j}$ is independent on the choice of
$i,j$ and lies in the kernel of $SK_1(R,U;d)$. This element is called a
Steinberg symbol and is denoted by $\{\bar u,\bar v\}$.
By definition each Steinberg symbol can be written as a product of
$13$ elementary matrices.

For elements $u_i, u_i' \in R$, $i=1,\dots,p \leq d/2$ such that $u_iu_i' = 1 \mod U$
consider the elements
$$
\widetilde{w}(u_i;u'_i) = \prod_i e_{i,p+i}(u_i) . \prod_i e_{p+i,i}(-u'_i). \prod_i e_{i,p+i}(u_i)
$$
and
$$
\widetilde{h}(u_i;u'_i):=\widetilde{w}(u_i;u'_i)\widetilde{w}(-1,\dots;-1,\dots)
$$
Using the definition of the Steinberg symbols it can be seen that
$$
\prod_i \{\bar u_i,\bar v_i\} = \widetilde{h}(u_iv_i;u'_iv'_i) \widetilde{h}(u_i;u'_i)^{-1} \widetilde{h}(v_i;v'_i)^{-1}.
$$
Therefore any product of $p = \lfloor d/2 \rfloor$ Steinberg symbols
can be written as a product of $13$ generalized elementary matrix --
here we have used that $\prod_i e_{i,p+i}(u_i)$ is a GEM.
This gives us that any element in $SK_1(R,U)$
can be written as
$$
13 \lceil (k+2) / p \rceil \leq 13(3+ k / p)  \leq 13 (3+ 3k/d)
$$
generalized elementary matrices.

Combining lemmas~\ref{congimage} and~\ref{kernel} we obtain:
\begin{theorem}
\label{gemcomm}
Let $H$ be a normal subgroup of finite index in $G_{d,k}$. Then
every matrix in $G_{d,k}=\SL_d(\Z[x_1,\dots,x_k])$ can be written as a product of
$77+39k/d$ generalized elementary matrices modulo $H$.
\end{theorem}

Finally we can complete the proof of Theorem~\ref{taucom}. Let $g \in G_{d,k}$,
using theorem~\ref{gemcomm} we can write $g$ as a product of $77+39k/d$ GEMs $v_s$
and an element in $h \in H = \ker \rho$. As in the proof of Theorem~\ref{t1} this
implies that
$$
||\rho(g)v-v|| \leq ||\rho(h)v-v|| + \sum_{s=1}^{77+39k/d} ||\rho(v_s)v-v|| <
                2(77 +39k/d)M_d(k)\epsilon.
$$

If $T_{d,k}= [\sqrt{2}(77+39k/d)M_d(k)]^{-1}$ then we will have that $||\rho(g)v-v|| < \sqrt{2}$ for
any $g\in G$ and the representation $\rho$ will have an invariant vector,  which completes the
proof that the group $\SL_d(R)$ has the property \emph{T} and also provides a
lower bound for the \TauC:
$$
\TC(G_{d,k},S_{d,k}) \geq \frac{1}{\sqrt{2}(77+39k/d)M_d(k)}.
$$
It is easy to see that we have the inequality
$$
\sqrt{2}(77+39k/d)M_d(k) < 800\sqrt{d}\left(1+(k/d)^{3/2}\right),
$$
which completes the proof of theorem~\ref{taucom}.$\square$

\medskip

Theorems~\ref{t1} and~\ref{gemcomm} are now proved Lemmas~\ref{abe} and~\ref{congimage}.

\section{Relative property \emph{T} of $\EL_2(R) \ltimes {R}^2,{R}^2$}
\label{relative}

In this section we shall
show that the pair of groups $(\EL_2(R) \ltimes {R}^2, {R}^2)$ has
relative property \emph{T}, and we shall estimate its relative Kazhdan
constant. This result was first proven by M.~Burger~\cite{Bur}
in the case $R=\mathbb{Z}$ and later generalized by Y.~Shalom~\cite{YSh} to
the case of finitely generated commutative rings. In their proofs they use that
$\SL_2(\Z)$ is generated by the elementary matrices (and is thus
equal to $\EL_2(\Z)$), but this assumption is not necessary. The argument
which follows gives a better estimate of the relative
Kazhdan constant in the case of commutative ring $R$ and
is based on ideas from~\cite{K}.

First we start with a definition of the \emph{relative property T} for finitely generated
discrete groups:

\begin{definition} A pair of groups $(G,H)$ such that $H\leq G$, where $G$ is
generated by a finite set $S$  is said to
have the \emph{Relative Kazhdan Property T} if there is
$\epsilon  >0$
such that for every irreducible (continuous) unitary representation
$\rho: \ G \rightarrow U(\mr{H})$ on a
Hilbert space $\mr{H}$ without $H$ invariant vectors and
for every  vector $v \not = 0$ in $\mr H$ there is some
$s \in S$ such that $||\rho(s)v-v|| > \epsilon ||v||$.
The largest $\epsilon$ with this property is called the
\emph{Relative \KaC} and is denoted by $\mr{K}(G,H;S)$. It is
easy to see that the relative property \emph{T} depends only on the groups
$G$ and $H$ and is independent of the generating set $S$, however the
value of the relative Kazhdan constant depends also on $S$.
\end{definition}

We will show that $\left(\EL_2(R^n_k) \ltimes {R^n_k}^2, {R_k^n}^2\right)$ has the relative property \emph{T}
by proving the following theorem (which is inspired by~\cite{YSh}):

\begin{theorem}
\label{relT}
Let $F_1$ be set of $4(k+1)$ elementary matrices in $\EL(2;R_k^n)$ with $\pm 1$ and $\pm x_j$ off the
diagonal and $F_2$ be the set of standard basis vectors in ${R_k^n}^2$ and their inverses.
Let $\rho: \Gamma \rightarrow U(\mr{H})$ be a unitary representation of
the group
$\Gamma:=\EL_2(R_k^n) \ltimes {R_k^n}^2$ in a Hilbert space $\mr{H}$.
Let $v\in \mr{H}$ be a unit vector
such that $||\rho(g)v-v|| <\epsilon=M(k)^{-1}$ for all $g\in F_1\cup F_2$,
then $\mr{H}$ contains an ${R^n_k}^2$-invariant vector.
Here the numbers $M(k)$ are defined by
$M(0) = 2+\sqrt{10}$ and
$$
\begin{array}{r@{}c@{}l}
M(k)^2 &{}={}&
\displaystyle
4M(0)^2 +  (6\sqrt{k}+4)\sqrt{(M(0)+\sqrt{k})^2+k+2\sqrt{k}} + \\
& & \displaystyle \quad + 6 \sqrt{k}M(0) + 7k+10\sqrt{k}+2.
\end{array}
$$
This shows that
$\EL_2(R_k^n) \ltimes {R_k^n}^2, {R_k^n}^2$ has relative property \emph{T} and
$$
\KC(\EL_2(R_k) \ltimes R_k^2, R_k^2;F_1\cup F_2) \leq 1/{M(k)}.
$$
Using the definition of $M(k)$ it can be verified that $M(0) < 5.17$, $M(1)<14.92$, $M(2)< 16.47$ and
$M(k) \leq \sqrt{18}(\sqrt{k} + 3)$.
\end{theorem}

Before starting the proof of this theorem, let us note that there
are two natural actions of $\EL_2(R)$ on $R^2$, a `left' action where we think of
the elements in $R^2$ as column vectors, i.e., we embed $\EL_2(R) \ltimes R^2$
in $\EL_3(R)$ as follows
$$
(A,B) \to \left( \begin{array}{cc} A & B \\ 0 & 1 \end{array} \right).
$$
There is also a `right' one where we think of $R^2$ as a row vectors and we use
the embedding
$$
(A,B) \to \left( \begin{array}{cc} A & 0 \\ B & 1 \end{array} \right).
$$
In the case of commutative ring these two actions give isomorphic semidirect products, which is
not always the case if the ring $R$ is not commutative. However we have that the `left'
semi-direct product is isomorphic to the `right' semidirect product over the opposite ring $R^{opp}$.
In the proof we will consider only the `left' action. By the previous remark
the theorem is valid also for the `right' semidirect products
and the proof is essentially the same.

\noindent
\textbf{Proof of Theorem~\ref{relT}:}
Assume that we have a representation $\rho$ which satisfies all the conditions in the theorem and
the Hilbert space $\mr{H}$ does not contain
an $R^2$ invariant vector.
Let $P$ be the 
projection valued measure on the dual $(R^2)^* = \Rd$, coming
form the restriction of the representation $\rho$ to $R^2$.
Here by $R^*$ we have denoted the dual of the additive group of the ring $R$, i.e.,
$
R^* = \Hom \left(R+, S^1\right),
$
where $S^1=\{z\in \C \mid |z|=1\}=U(\C)$ is the unit circle in the complex plane.
The set $R^*$ is an abelian group with natural topology.
For $t\in R^*$ and $r\in R$ by $t(r)$ we will denote the natural pairing of $R^*$ and $R$.
Let $\mu_v$ be the probability measure on $\Rd$, defined by %
$\mu_v(B)=\la P(B)v,v \ra$. The measure
$\mu_v$ is supported on
${R^*}^2\setminus \{(0,0)\}$, because by assumption $\mr{H}$ does not
contain an $R^2$ invariant vector and by the construction
of the measure $P$, $P(\{0\})$ is the
projection onto the subspace of $R^2$ invariant vectors in $\mr{H}$.

For an element $T \in \Rd$ we will write $T=(t_1,t_2)^t$,
where  $t_1$ and $t_2$ are in $R^*$.

\begin{lemma}
\label{inv}
Let $P_i = \{T \in \Rd \mid \,\, \Re(t_i(1)) < 0 \}$, then %
$\mu_v(P_{i}) \geq 1-\epsilon^2/2$.
\end{lemma}
\textbf{Proof:}
By the definition of the measure $\mu_v$, we have
$$
|| \rho(g_{i})v-v||^2 =
\int_{\Rd} | t_i(1) -1 |^2\, d\mu_v,
$$
where $g_1,g_2\in G$ form the standard basis of $R^2$. By assumption $v$ is
almost preserved by $\rho(g_{i})$ therefore $|| \rho(g_{i})v-v||^2 \leq \epsilon$.
If we break the integral into two integrals over $P_i$ and its complement, we get
$$
\epsilon^2 \geq \int_{\Rd\setminus P_i} \!\!\!\!\!| t_i(1) -1 |^2\, d\mu_v +
\int_{P_i} | t_i(1) -1 |^2\, d\mu_v
\geq \int_{\Rd\setminus P_i} \!\!\!\!\! 2 \, d\mu_v = 2(1-\mu_v(P_i)),
$$
which gives that the measure of the set $P_i$ is large enough. $\square$

\begin{lemma}
\label{invariance}
For every measurable set $B \subset \Rd$ and
every elementary matrix $g \in F_1$, we have that
$$
|\mu_v(gB) - \mu_v(B)| \leq 2\epsilon\sqrt{\mu_v(B)}+ \epsilon^2
\quad \mbox{or} \quad
|\sqrt{\mu_v(gB)} - \sqrt{\mu_v(B)}| \leq \epsilon.
$$
The action of $\EL_2(R)$ on $\Rd$ is the dual of the action of $\EL_2(R)$ on
$R^2$.
\end{lemma}

\textbf{Proof:}
Using the properties of the projection valued measure $P$, we have
$$
\begin{array}{r@{\,\,}l}
 |\mu_v(gB) - \mu_v(B)| = & 
       | \la\pi(g^{-1})P(B)\pi(g)v,v\ra - \la P(B)v,v\ra| \leq \\
 \leq & \displaystyle \rule[15pt]{0pt}{0pt}
       |\la\pi(g^{-1})P(B)(\pi(g)v-v),v\ra | +
       |\la P(B)v,(\pi(g)v-v\ra)| = \\
 = & \displaystyle \rule[15pt]{0pt}{0pt}
       2|\la\pi(g)v-v,P(B)v\ra| +
       \la P(B)(\pi(g)v-v),\pi(g)v-v\ra \leq  \\
 \leq & \displaystyle \rule[15pt]{0pt}{0pt}
       2\epsilon \sqrt{ \mu_v(B)} + \epsilon^2,
\end{array}
$$
where the final inequality follows from the facts
that $v$ is $(F,\epsilon)$ invariant
vector and $||P(B)v||^2 = \mu_v(B)$. $\square$

We will need the following slight generalization of the above lemma:
\begin{lemma}
\label{minvariance}
Let $A$ and $B$ are measurable sets in $\Rd$. Suppose that $A$ decomposes
as a disjoint union of the sets $A_i$ for $i=1,\dots, s$ and there
exist elements $g_i \in F_1$ such that the sets $B_i=g_i(A_i)$ are disjoint
subsets of $B$ then
$$
\sqrt{\mu_v(A)} \leq \sqrt{\mu_v(B)} + \sqrt{s}\epsilon.
$$
\end{lemma}

\textbf{Proof:}
Applying lemma~\ref{invariance} to the sets $A_i$ yields to
$$
\begin{array}{r@{}c@{}l}
\mu_v(A) & {}={} & \displaystyle
\sum_{i=1}^s \mu_v(A_i) \leq \sum_{i=1}^s \left[ \mu_v(B_i) + 2\epsilon \sqrt{ \mu_v(B_i)} + \epsilon^2 \right] \leq \\
& \leq & \displaystyle
\sum\mu_v(B_i) + 2\epsilon\sqrt{s\sum \mu_v(B_i)} + s\epsilon^2 \leq \\
& \leq & \displaystyle
\mu_v(B) + 2\epsilon\sqrt{s\mu_v(B)} + s\epsilon^2 = \left(\sqrt{\mu_v(B)} + \sqrt{s}\epsilon\right)^2,
\end{array}
$$
where we have used that for nonnegative numbers $a_i$ we have the inequality
$$
\displaystyle \sqrt {\sum_{i=1}^s a_i} \leq \sum_{i=1}^s \sqrt{a_i} \leq \sqrt {s \sum_{i=1}^s a_i}. \square
$$

\medskip

Let us first consider the case when $k=0$ and $R=\Z$.
In this case the dual
of the ring $R$ is just the unit circle in complex plane. Also we can
think of $\Rd$ as the torus $\T^2 = S^1 \times S^1 = \R^2/\Z^2$ which we will identify
with the square $(-1/2,1/2] \times (-1/2,1/2]$. It is well known that there are only a few
$\SL_2(\Z)$ invariant measures on the torus $\T^2$. The next lemma is a quantitative
version of this fact.
\begin{lemma}
\label{R2}
Let $\mu$ be a finitely additive measure on $\T^2$ such that:
\begin{enumerate}
\item %
$\mu(|x| \geq 1/4) \leq \epsilon^2/2$ and %
$\mu(|y| \geq 1/4) \leq \epsilon^2/2$, %
\item %
$|\mu(gB) - \mu(B)| \leq 2\epsilon \sqrt{\mu(B)} + \epsilon^2$ %
for any Borel set $B$ and any elementary matrix $g\in F_1\cap \SL_2(\Z)$.
\end{enumerate}
Then the measure $\mu$ satisfies:
$$
\mu(\T^2\setminus \{(0,0)\}) \leq (2+\sqrt{10})^2\epsilon^2.
$$
\end{lemma}
\textbf{Proof:}
Let us define the Borel subsets $A_i$ and $A_i'$ of $\T^2$ using
the picture:

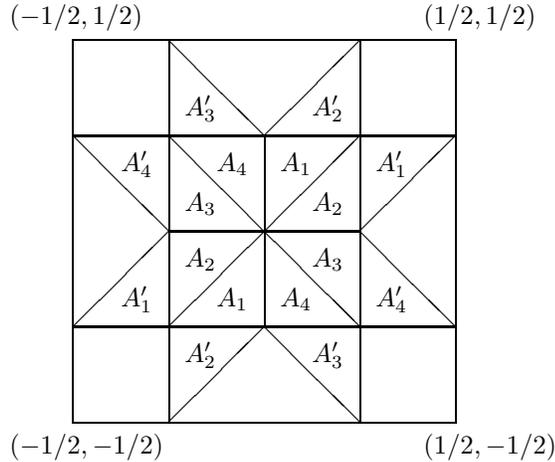
\setlength{\unitlength}{3947sp}%
\begin{picture}(3000,2900)(200,-2350)
\thinlines
\put(1800,-2100){\framebox(2400,2400){}}
\put(1800,-1500){\line(1,0){2400}}
\put(1800, -300){\line(1,0){2400}}
\put(2400, -900){\line(1,0){1200}}
\put(2400,-2100){\line(0,1){2400}}
\put(3600,-2100){\line(0,1){2400}}
\put(3000,-1500){\line(0,1){1200}}
\put(2400,-1500){\line(1, 1){1200}}
\put(2400, -300){\line(1,-1){1200}}
\put(1800,-1500){\line(1, 1){600}}
\put(1800, -300){\line(1,-1){600}}
\put(3600, -900){\line(1, 1){600}}
\put(3600, -900){\line(1,-1){600}}
\put(2400,-2100){\line(1, 1){600}}
\put(2400,  300){\line(1,-1){600}}
\put(3000,-1500){\line(1,-1){600}}
\put(3000, -300){\line(1,1){600}}
\put(3100,-520){$A_1$}
\put(3300,-770){$A_2$}
\put(3300,-1120){$A_3$}
\put(3100,-1370){$A_4$}
\put(2700,-1370){$A_1$}
\put(2500,-1120){$A_2$}
\put(2500,-770){$A_3$}
\put(2700,-520){$A_4$}
\put(3700,-520){$A_1'$}
\put(3300,-170){$A_2'$}
\put(3300,-1720){$A_3'$}
\put(3700,-1370){$A_4'$}
\put(2100,-1370){$A_1'$}
\put(2500,-1720){$A_2'$}
\put(2500,-170){$A_3'$}
\put(2100,-520){$A_4'$}
\put(1400,-2300){$(-1/2,-1/2)$}
\put(1400, 400){$(-1/2,1/2)$}
\put(4000,-2300){$(1/2,-1/2)$}
\put(4000,400){$(1/2,1/2)$}
\end{picture}

\noindent
Each set $A_i$ or $A_i'$ consists of the interiors of two
triangles and part of their boundary (not including the vertices).
The sets $A_i$ do not contain the side which is part of the small
square, they also do not contain their clockwise boundary
but contain the counter-clockwise one. Each set $A_i'$ includes
only the part of its boundary which lies on the small square.

It can be seen, from the picture, that the elementary matrices
$e_{ij}^{\pm} = e_{ij}(\pm 1)\in F_1$,
act on the sets $A_i$ as follows:
$$
\begin{array}{ll}
e_{12}^+(A_3 \cup A_4')  = A_3 \cup A_4 \quad &
e_{21}^+(A_3'\cup A_4 )  = A_3 \cup A_4 \\
e_{12}^-(A_1'\cup A_2 )  = A_1 \cup A_2  \quad &
e_{21}^-(A_1 \cup A_2')  = A_1 \cup A_2.
\end{array}
$$
Using the properties of the measure $\mu$ the above equalities
imply the inequalities:
$$
\begin{array}{l}
\mu(A_1) + \mu(A_2) \leq
     \mu(A_1') + \mu(A_2) + \epsilon^2 +
     2\epsilon\sqrt{\mu(A_1') + \mu(A_2 )} \\
\mu(A_1) + \mu(A_2) \leq
     \mu(A_1) + \mu(A_2') + \epsilon^2 +
     2\epsilon\sqrt{\mu(A_1 ) + \mu(A_2')} \\
\mu(A_3) + \mu(A_4) \leq
     \mu(A_3') + \mu(A_4) + \epsilon^2 +
     2\epsilon\sqrt{\mu(A_3') + \mu(A_4 )} \\
\mu(A_3) + \mu(A_4) \leq
     \mu(A_3) + \mu(A_4') + \epsilon^2 +
     2\epsilon\sqrt{\mu(A_3 ) + \mu(A_4')}. \\
\end{array}
$$
Adding these inequalities and noticing that
$$
\begin{array}{l}
\mu(A_1') + \mu(A_4') \leq
     \mu(\{ |x| \geq 1/4\}) \leq
     \epsilon^2/2 \,\,\, \mbox{ and}\\
\mu(A_2') + \mu(A_3') \leq
     \mu(\{ |y| \geq 1/4\}) \leq
     \epsilon^2/2
\end{array}
$$
we obtain
$$
\sum \mu(A_i) \leq
    4\epsilon^2 + \sum \mu(A_i') +
        2\epsilon \sqrt{4\left(\sum \mu(A_i) +
    \sum \mu(A_i')\right)} \leq
$$
$$
    \leq 5\epsilon^2 + 4\epsilon \sqrt{\sum \mu(A_i) + \epsilon^2}.
$$

After substituting $c =\sqrt{\sum \mu(A_i) + \epsilon^2}$
and solving the resulting quadratic inequality one obtains
$\sum\mu(A_i) \leq (13 +4 \sqrt{10}) \epsilon^2$.
Finally we can use that
$$
\mu(\T^2\setminus \{(0,0)\}) \leq \sum \mu(A_i) + \mu(\{ |x| \geq
1/4\}) + \mu(\{ |y| \geq 1/4\}) \leq
$$
$$
\leq (14 + 4 \sqrt{10}) \epsilon^2 = (2+\sqrt{10})^2\epsilon^2,
$$
which completes the proof of the lemma. $\square$

This completes the proof of Theorem~\ref{relT} in the case $k=0$ ---
by assumption the Hilbert space $\mr{H}$ does not contains
$R^2$ invariant vectors. However, Lemma~\ref{R2} gives us that
$\mu(0,0)>0$ if $\epsilon< \frac{1}{2+\sqrt{10}}=M(0)^{-1}$ which is a contradiction
with the assumption the $\mr{H}$ does not contain $\Z^2$ invariant vectors.

\medskip

In order to prove the general case of Theorem~\ref{relT}
we need to generalize Lemma~\ref{R2} to the case of larger rings.
Recall, that $R_k^n$ denotes the free ring $\Z\la x_1,\dots,x_k \ra$. Let $F$ denote
the set of elementary
matrices in $\EL_2(R_k^n)$, with $\pm 1$ and $\pm x_i$ off the diagonal.
Denote the dual of the abelian group $R_k^n= \Z\la x_1,\dots,x_k\ra$ by ${R_k^n}^*$.
\begin{lemma}
\label{R2k}
Let $\mu$ be a finitely additive measure 
on $\left.{R^n_k}^*\right.^2$,
such that
\begin{enumerate}%
\item %
$\mu(\{T\mid \Re(t_1(1)) < 0 \}) \leq \epsilon^2/2$ and %
$\mu(\{T\mid \Re(t_2(1)) < 0 \}) \leq \epsilon^2/2$. %
Here $t_i$ denote the components of the element $T\in \left.{R^n_k}^*\right.^2$.

\item %
$|\mu(gB) - \mu(B)| \leq 2\epsilon \sqrt{\mu(B)} + \epsilon^2$ %
for any Borel set $B\subset \left.{R^n_k}^*\right.^2$ and any elementary matrix %
$g\in F_1 \subset \EL_2(R_k)$.
\end{enumerate}
Then 
we have
$$
\mu(\left.{R^n_k}^*\right.^2\setminus \{(0,0)\}) \leq %
M(k)^2 \epsilon^2, 
$$
where the number $M(k)$ is defined in Theorem~\ref{relT}.
\end{lemma}
\textbf{Proof:}
The idea for the proof is based on the proof of Lemma 3.3
from~\cite{YSh} - the main difference is that we do
not use induction on $k$ and do everything in one step, which allows us
to drop the assumption that the ring $R$ is commutative.

By the Lemma~\ref{R2} we know that if we only consider the action of
$\Z\subset R^n_{k}$ we have that
$$
\mu(\left.\Z^*\right.^2\setminus \{(0,0)\}) \leq M(0)^2\epsilon^2.
$$
There is a natural restriction map $\mathrm{res}: {R^n_{k}}^* \to \Z^*$ coming from
the inclusion $\Z\subset R_k^n$ and the correct way of writing the
above inequality is
$$
\mu(\left.{R^n_k}^*\right.^2 \setminus \mathrm{res}^{-1}\{(0,0)\}) \leq M(0)^2\epsilon^2.
$$

By definition, we have that $R^n_{k}= \Z\la x_1,\dots,x_{k}\ra$, which
allows us to write
$$
{R^n_{k}}^* \simeq \Z^*\la\la x_1^{-1},\dots,x_k^{-1}\ra\ra
$$
as a topological space. Here $\Z^*\la\la x_1^{-1},\dots,x_k^{-1}\ra\ra$
is the formal power series on non commuting variables $x_k^{-1}$ with coefficients in $\Z^*$.
This isomorphism gives a valuation $\nu$ on ${R^n_{k}}^*$,
corresponding to the total degree of the
leading term.

Let us define the following subsets of $\left.{R^n_{k}}^*\right.^2\setminus\{0,0\}$:
$$
\begin{array}{r@{}c@{}l}
A & {} = {} & \left\{ (\chi_1,\chi_2) \mid \nu(\chi_1) > \nu(\chi_2) > 0 \right\} \\
B &    =    & \left\{ (\chi_1,\chi_2) \mid \nu(\chi_1) = \nu(\chi_2) > 0 \right\} \\
C &    =    & \left\{ (\chi_1,\chi_2) \mid \nu(\chi_2) > \nu(\chi_1) > 0 \right\} \\
D &    =    & \left\{ (\chi_1,\chi_2) \mid \nu(\chi_1)=0 \mbox{ or } \nu(\chi_2) = 0 \right\} \\
\end{array}
$$
Lets consider the action of $e_{i,j}(r), r\in R_k^n$ on the element $(\chi_1,\chi_2)\in \left.{R^n_{k}}^*\right.^2$,
by the definition of this action we have
$$
\left(e_{1,2}(r)\circ(\chi_1,\chi_2) \right)\left(\begin{array}{@{}c@{}}r_1 \\r_2\end{array}\right) =
(\chi_1,\chi_2)\left(e_{1,2}(r)^{-1} \left(\begin{array}{@{}c@{}}r_1 \\r_2\end{array}\right)\right) =
(\chi_1,\chi_2)\left(\begin{array}{@{}c@{}}r_1 -r\cdot r_2\\r_2\end{array}\right)
$$
and
$$
\left(e_{2,1}(r)\circ(\chi_1,\chi_2) \right)\left(\begin{array}{@{}c@{}}r_1 \\r_2\end{array}\right) =
(\chi_1,\chi_2)\left(e_{2,1}(r)^{-1} \left(\begin{array}{@{}c@{}}r_1 \\r_2\end{array}\right)\right) =
(\chi_1,\chi_2)\left(\begin{array}{@{}c@{}}r_1 \\r_2-r\cdot r_1\end{array}\right)
$$
This shows that
$$
e_{1,2}(1)\left(A\right) \subset B  \quad
e_{2,1}(1)\left(C\right) \subset B.
$$
and
the first two inequalities can be written as
$$
\sqrt{\mu(A)} \leq \sqrt{\mu(B)} + \epsilon
\quad \mbox{and} \quad
\sqrt{\mu(C)} \leq \sqrt{\mu(B)} + \epsilon.
$$

\noindent
\textbf{Claim:} The union $A\cup B$ can be decomposed as a disjoint union of
$k$ sets $A_i$ such that the sets $\widetilde{A_i}= e_{1,2}(x_i)\left(A_i\right)$ are disjoint and
lie in $C\cup D$. Similarly, $C\cup B$ can be written as $\cup C_i$ such that
$\widetilde{C_i}= e_{2,1}(x_i)\left(C_i\right)$ are disjoint subsets of $A\cup D$.

\noindent
\textbf{Proof of the Claim:}
Let $(\chi_1,\chi_2)\in A\cup B$.
Let $\bar\chi$ be the leading term of $\chi_2$. Then $\bar\chi$ is a non trivial element of the
dual of the space of all homogeneous polynomials of degree $\nu(\chi_2)>0$ in $x_1,\dots,x_k$.
Let $\bar\chi_{(i)}$ be the element in the dual of the space of all homogeneous polynomials of
degree $\nu(\chi_2) -1 $ defined by
$$
\bar\chi_{(i)}(f) = \bar\chi (x_i f).
$$
Let $s(\bar\chi)$ be the smallest index such that $\bar\chi_{(i)}$ is not trivial -- such index exists
because $\bar\chi$ is non zero. Define the sets $A_i$ by
$$
A_i = \left\{ (\chi_1,\chi_2)\in A\cup B \mid s(\bar\chi) = i \right\}.
$$
It is clear that $A \cup B$ is a disjoint union of $A_i$ and that
$\widetilde{A_i}= e_{1,2}(x_i)\left(A_i\right)$
are disjoint subsets of $C\cup D$.

The second part of the clam is proves is the same way
but one uses the first component
$\chi_1$ of $\bar\chi$ not the second one $\chi_2$.
$\square$

We are ready to complete the proof of Lemma~\ref{R2k}.
Applying Lemma~\ref{minvariance} for the sets $A\cup B$ and $C\cup B$ gives us:
$$
\begin{array}{r@{}c@{}l}
\sqrt{\mu(A\cup B)} &{}\leq{}& \sqrt{\mu(C \cup D)} + \sqrt{k}\epsilon\\
\sqrt{\mu(C\cup B)} &{}\leq{}& \sqrt{\mu(A \cup D)} + \sqrt{k}\epsilon.
\end{array}
$$
Squaring and adding the the above inequalities leads to
$$
\begin{array}{r@{}c@{}l}
2\mu(B) & {} \leq {} & \displaystyle \rule[15pt]{0pt}{0pt}
2\mu(D) + 2 \epsilon\sqrt{k(\mu(A) + \mu(D))} + 2 \epsilon\sqrt{k(\mu(A) + \mu(D))} + 2k\epsilon^2  \leq \\
&\leq & \displaystyle \rule[15pt]{0pt}{0pt}
2\mu(D) + 4 \epsilon\sqrt{k\mu(D)} + 2 \epsilon\sqrt{k\mu(A)} +2 \epsilon\sqrt{k\mu(C)} + 2k\epsilon^2 \leq \\
&\leq & \displaystyle \rule[15pt]{0pt}{0pt}
2\mu(D) + 4 \epsilon\sqrt{k\mu(D)} + 4 \epsilon\sqrt{k\mu(B)} + (2k+4\sqrt{k})\epsilon^2 \\
\end{array}
$$
The last inequality can be rewritten as
$$
(\sqrt{\mu(B)} - \sqrt{k}\epsilon)^2 \leq (\sqrt{\mu(D)}+\sqrt{k}\epsilon)^2 + (k+2\sqrt{k})\epsilon^2
$$
Therefore we have
$$
\begin{array}{c@{}c@{}l}
\sqrt{\mu(A)} & {} \leq {} & \left(\sqrt{(M(0)+\sqrt{k})^2+ (k+2\sqrt{k})} + \sqrt{k}+1\right) \epsilon \\
\sqrt{\mu(B)} &  \leq      & \left(\sqrt{(M(0)+\sqrt{k})^2+ (k+2\sqrt{k})} + \sqrt{k} \right) \epsilon\\
\sqrt{\mu(C)} &  \leq      & \left(\sqrt{(M(0)+\sqrt{k})^2+ (k+2\sqrt{k})} + \sqrt{k}+1\right) \epsilon\\
\end{array}
$$
Finally we have
$$
\begin{array}{r@{}c@{}l}
\displaystyle
\mu(\left.R^*_{k}\right.^2\setminus \{(0,0)\})&{} = {}&
\displaystyle
\mu(A \cup B \cup C \cup D) \leq \mu(A) + \mu(B) + \mu(C) + \mu(D) \leq \\
& \leq & \displaystyle
\left( 4M(0)^2 + (6\sqrt{k}+4)\sqrt{(M(0)+\sqrt{k})^2+ k+2\sqrt{k}} +  \right. \\
&& \displaystyle \quad
 \left. \rule[15pt]{0pt}{0pt}6 \sqrt{k}M(0) +  7k+10\sqrt{k}+2\right) \epsilon^2 =
 M(k)^2 \epsilon^2,
\end{array}
$$
which completes the proof of lemma~\ref{R2k}. $\square$

\bigskip

Finally Theorem~\ref{relT} follows easily:
by assumption the Hilbert space $\mr{H}$ does not contain
$R^2$ invariant vectors. However lemma~\ref{R2k} gives us that
$\mu(0,0)>0$ if $\epsilon< M(k)^{-1}$ which is a contradiction.

Using that $M(0) = 2+\sqrt{10}$, it can be easily shown that
$$
M(0) < 5.17, \quad M(1) < 14.92,  \quad M(2) < 16.47, \quad\mbox{and} \quad M(k) < 3\sqrt{2}(\sqrt{k}+3).
\square
$$

\begin{corollary}
\label{abelian1}
Let $\rho:\Gamma \to U(\mr{H})$ be a unitary representation of the group
$$
\Gamma= \EL_2(R^n_k)\ltimes {R^n_k}^2.
$$
Let $v\in \mr{H}$ be a unit vector such that
such that $||\rho(g)v-v|| <\epsilon$ for $g\in F_1\cup F_2$.
Then for every $g$ in ${R^n_k}^2$
we have $||\rho(g)v - v || \leq 2M(k)\epsilon$.
\end{corollary}
\textbf{Proof:}
Let us decompose the Hilbert space $\mr{H}$ as the orthogonal direct sum $\mr{H}_0 \oplus \mr{H}_1$, where $\mr{H}_0$
contains all ${R^n_k}^{2}$ invariant vectors and $\mr{H}_1=\mr{H}_0^\perp$. We have that both $\mr{H}_0$
and $\mr{H}_1$ are closed under the action of the group $\Gamma$,
because ${R^n_k}^2$ is a normal subgroup of $\Gamma$.

Let us write
$v=v_0+v_1$, where $v_i \in \mr{H}_i$. Since there are no ${R^n_k}^2$-invariant vectors in $\mr{H}_1$,
there exists $h \in F \cup G$ such that %
$||\rho(h)v_1 - v_1|| \geq M(k)^{-1}||v_1||$. But, we have that
$$
||\rho(h)v - v||^2 = ||\rho(h)v_0 - v_0||^2 + ||\rho(h)v_1 - v_1||^2
    \leq \epsilon^2.
$$
The final inequality implies $||v_1|| \leq M(k) \epsilon$. For any $g\in {R^n_k}^2$, we
have
$$
||\rho(g)v - v||^2 = ||\rho(g)v_0 - v_0||^2 + ||\rho(g)v_1 - v_1||^2
   \leq 0 + 4 || v_1||^2 \leq 4M(k)^2\epsilon^2,
$$
therefore $||\rho(g)v - v || \leq 2M(k)\epsilon$. $\square$

As another corollary we obtain part a) of Lemma \ref{abe} from Section \ref{outline}.

\bigskip

Using the similar methods we can estimate the relative Kazhdan constant for the pairs
$$
\EL_p(R) \ltimes R^{p}, R^{p},
$$
with $p\geq 2$ and
$$
\left( \EL_p(R) \times \EL_q(R) \right) \ltimes R^{pq}, R^{pq},
$$
where $p\geq 2$. Here we consider the elements of $\left( \EL_p(R) \times \EL_q(R) \right) \ltimes R^{pq}$ as
$(p+q) \times (p+q)$ matrices over $R$ of type
$$
\left(
\begin{array}{cc} A & C \\ 0 & B \end{array}
\right),
$$
where $A\in \EL_p(R)$, $B \in \EL_q(R)$ and $C \in R^{pq}$.
Equivalent way of say this is the we consider
$R^{pq}$ as a left $\EL_p(R)$-module and a right $\EL_q(R)$-module.

\begin{theorem}
\label{relTpq}
a)
Let $F_{p;k}$ be the following generating set of
$\Gamma_p\!=\EL_p(R_k^n) \! \ltimes \!\left.R_k^n\right.^{p}$
$$
\begin{array}{c@{}c@{}l}
F_{p;k}  & {}={} &
\left\{\Id \pm e_{ij} \in \EL_p(R) \ | \ 1\leq i\not=j \leq p \right\}
\cup \\
& & \quad \cup
\left\{\Id \pm x_l\cdot e_{ij} \in \EL_p(R)\ | \ |i-j|=1, 1\leq l\leq k \right\}
\cup \\
& &\quad  \cup
\left\{\pm 1_{i} \in {R}^{p}\ | \ 1\leq i\leq p \right\}
\end{array}
$$
consisting of some
elementary matrices and the standard generators of $R^{p}$.

Let $\rho: \Gamma_p \rightarrow U(\mr{H})$ be a unitary representation of $\Gamma_p$ acting
on a Hilbert space $\mr{H}$. Let $v\in \mr{H}$ be a unit vector
such that $||\rho(g)v-v|| <\epsilon=M(k;p)^{-1}$ for all $g\in F_p$,
then $\mr{H}$ contains an $R_k^p$-invariant vector.
Here $M(k;p)$ denotes the number
$$
M(k;p)= \sqrt{2M(k)^2 + 2M(k)\sqrt{p-2} + p-2},
$$
where $M(k)$ is the constant defined in Theorem~\ref{relT}. \bigskip

b)
Let $F_{p,q;k}$ be the generating set of
$\Gamma_{p,q}=\left(\EL_p(R_k^n) \times \EL_q(R_k^n)\right) \ltimes \left.R_k^n\right.^{pq}$:
consisting of the elementary matrices in $\EL_p(R)$ and $\EL_q(R)$ and the
standard generators
of $R^{pq}$.

Let $\rho: \Gamma_{p,q} \rightarrow U(\mr{H})$ be a unitary representation of $\Gamma_{p,q}$ on
a Hilbert space $\mr{H}$. Let $v\in \mr{H}$ be a unit vector
such that $||\rho(g)v-v|| <\epsilon=M(k;p,q)^{-1}$ for all $g\in F_{p,q}$,
then $\mr{H}$ contains an $H$-invariant vector.
Here $M(k;p,q)$ denotes the number
$$
M(k;p,q)= \sqrt{2M(k;p)^2 + 2M(k;p)\sqrt{q-1} + q-1}.
$$
\end{theorem}
\textbf{Proof:}
a) The proof is similar to to the one of Theorem~\ref{relT}.
Let $\mu_v$ be the measure on $\left.R^*_k\right.^p$ coming from the restriction of the
representation $\rho$ to the the abelian group $R_k^p$.
As in the case $p=2$ the measure $\mu_v$
is almost preserved under the action of $\EL_p(R)$.

We will write an element  $T \in \left.R^*_k\right.^{p}$ as $(t_1,\dots,t_p)^t$,
where $t_i \in R^*_k$.
Using the restriction $\left.R^*_k\right.^p \to\left.R^*_k\right.^2$ coming form the
inclusion $R_k^2 \subset R_k^p$ and Lemma~\ref{abelian1} we can see that
$$
\mu_v(\{T \mid t_1\not= 0 \mb{ or } t_2\not= 0 \}) \leq M(k) \epsilon^2.
$$

Let us define the Borel subsets $B_i,C_i$ of $\left.R^*_k\right.^p$ by
$$
\begin{array}{r@{}c@{}l}
B_i &{}= {}&\{ T \mid t_k =0 \mbox{ for } k\leq i\}, \mbox{ and}\\
C_i &= &\{ T \mid t_1 = t_i \not =0, t_k = 0 \mbox{ for } 1<k<i\}.
\end{array}
$$
The elementary matrix $g_{1i} \in \SL_p\subset \EL_p(R)$ sends $B_{i-1} \setminus
B_i$ into $C_i$ for any $i\geq 3$.
Now, notice that the sets $C_i$ are
disjoint for $i=2,\dots, p$, and their union lies in the set %
$$
C:=\{ T \mid t_1 \not =0, t_2 = 0 \} \subset \{ T \mid t_1\not= 0 \mb{ or } t_2\not= 0 \}.
$$
Therefore, by lemma~\ref{minvariance}
we have
$$
\sqrt{\mu_v(B_2 \setminus B_p )} \leq \sqrt{ \mu_v(C)} + \sqrt{(p-2)}\epsilon.
$$
We know that the measure of $C$ is very small which gives
$$
\begin{array}{r@{\,\,}l}
\mu_v(\left.R^*_k\right.^p\setminus & \{(0,\dots,0)\}) =
    \mu_v(\left.R^*_k\right.^p \setminus B_2) + \mu_v(B_2 \setminus B_p) \leq \\
\leq & \displaystyle \rule[15pt]{0pt}{0pt}
 M(k)^2\epsilon^2 + \left(M(k)\epsilon + \sqrt{p-2}\epsilon\right)^2
=  M(k;p)^2 \epsilon^2.
\end{array}
$$
By assumption $\epsilon < M(k;p)^{-1}$ which implies that the measure of the zero in
$\left. R_k^* \right.^p$ is
$\mu_v(\{(0,\dots,0)\}) >0$,
and this shows that $\mr{H}$ contains $R^p$-invariant vectors. \bigskip

b) The idea is the same as the one in part a):
Let $\mu_v$ be the measure on $\left.R^*_k\right.^{pq}$ coming from the restriction of the
representation $\rho$ to the abelian group $R_k^{pq}$. As in the case $q=1$ the measure $\mu_v$
is almost invariant under the action of $\EL_p(R)$ and $\EL_q(R)$.

We will write an element  $T \in \left.R^*_k\right.^{pq}$ as $(T_1,\dots,T_q)$ where
$T_i\in \left.R^*_k\right.^{p}$.
Using the restriction $\left.R^*_k\right.^{pq} \to\left.R^*_k\right.^p$ coming from the
inclusion $R_k^p \subset R_k^{pq}$ which is given by $T \to T_1$, and part a) we have
$$
\mu_v(\{T \mid T_1\not= 0) \leq M(k;p)^2 \epsilon^2 .
$$

Define the Borel subsets $B_i,C_i$ of $\left.R^*_k\right.^{pq}$ by
$$
\begin{array}{l}
B_i = \{ T \mid T_k =0 \mbox{ for } k\leq i\}, \mbox{ and}\\
C_i = \{ T \mid T_1 = T_i \not =0, T_k = 0 \mbox{ for } 1<k<i\}.
\end{array}
$$
The elementary matrix $g_{1i} \in \EL_q$ sends $B_{i-1} \setminus
B_i$ into $C_i$ for any $i$. Therefore, we have
$$
\mu_v(B_{i-1} \setminus B_i ) \leq
\mu_v(C_i) + \epsilon^2 + 2\epsilon\sqrt{\mu_v(C_i)}.
$$
The sets $C_i$, for $i=1,\dots, q$, are
disjoint and their union lies in the set %
$$
C=\{ T \mid T_1 \not =0 \}.
$$
Therefore, by lemma~\ref{minvariance} we have
$$
\sqrt{\mu_v(B_1 \setminus B_q )}  \leq
\sqrt{\mu_v(C)} + \sqrt{q-1}\epsilon.
$$
We know that the measure of $C$ is very small which gives us
$$
\begin{array}{r@{}c@{}l}
\mu_v(\left.R^*_k\right.^{pq}\setminus \{0\}) &{}={}& \displaystyle
    \mu_v(\left.R^*_k\right.^{pq} \setminus B_1) + \mu_v(B_1 \setminus B_q) \leq \\
& \leq & \displaystyle \rule[15pt]{0pt}{0pt}
 M(k;p)^2\epsilon^2 + \left(M(k;p)\epsilon + \sqrt{q-1}\epsilon\right)^2 = M(k;p,q)^2\epsilon^2\\
\end{array}
$$
By assumption $\epsilon < M(k;p,q)^{-1}$ which implies that $\mu_v(\{0\}) >0$.
This shows that $\mr{H}$ contains $R^{pq}$-invariant vectors. $\square$ \bigskip

\begin{corollary}
\label{abepq}
Let $\rho: G \to U(\mr{H})$ be a unitary representation of the group
$G=\EL(d;R)$
Let $v\in \mr{H}$ be a unit vector such that $||\rho(g)v-v|| <\epsilon$ for $g\in F_1\cup F_2$,
where
$F_1$ is the set of $2(d^2-d)$ elementary matrices with $\pm 1$ off the diagonal
and $F_2$ is the set $4k(d-1)$ elementary matrices $\Id \pm x_se_{ij}$ with $|i-j|=1$.
Then for every generalized elementary matrix $g$
we have $||\rho(g)v - v || \leq 2M_d(k)\epsilon$, where $M_d(k)=2M(k) + \sqrt{3d} \leq 6\sqrt{2}(\sqrt{k}+3) + \sqrt{3d}$.
\end{corollary}
\textbf{Proof:}
Use the inequality
$$
M(k;p,q) \leq 2M(k) + \sqrt{2p} + \sqrt{q} \leq 2M(k) + \sqrt{3d} = M_d(k)\mbox{ if }p+q\leq d. \square
$$

\medskip

\noindent
This corollary immediately implies part b) of Lemma~\ref{abe} from Section~\ref{outline}.

\section{Uniform bounded generation of $\SL_d(\bar{R})$}
\label{uni}

Lemma~\ref{congimage} is a well-known result for fields. In the case of
finite commutative rings the same proof works. Almost everything works
also in the case of ring satisfying Bass stable range condition.

\begin{definition}
A ring $R$ is said to satisfy the Bass stable range condition if for any
elements $a,b \in R$ such that $Ra + Rb = R$, there exist $t\in R$ such that
$R(a+tb)=R$.
\end{definition}
There is an equivalent definition:
\emph{If $a$ be an element $R$ and $I$ be a left ideal $R$, such that the
left ideal generated by $a$ and $I$ coincides with $R$. Then the there exist $i \in I$
such that the left ideal generated by $u=a+i$ is the whole ring $R$.
}

This is the left Bass stable range condition, there is also
a right version and it is known that the two conditions are equivalent.
Any finite ring satisfies the Bass stable range condition see~\cite{Va}.

If a ring $R$ satisfies the Bass stable range condition, there is an
algorithm which writes an element in $\EL_d(R)$ as a product of small
number of elementary matrices and an element in $\EL_{d-1}(R)$.

\begin{lemma}
\label{congimagest}
Let the ring $\bar R$ satisfied the Bass stable range condition.
Then any element $g\in \EL_d(\bar R)$ can be written as
$$
g= \prod_{s=1}^{d^2-1}E_s \cdot g_{11} \cdot  \prod_{t=1}^{(d^2-d)/2}E'_t
$$
where each $E_s, E'_t \in E(\bar R)$ is an elementary matrix and $g_{11}$ lies in the copy of
$\GL_1(\bar R)$ embedded in the upper right corner, i.e.,
$
g_{11} = e_{11}(a),
$
where $a$ is a unit in $\bar R$.
\end{lemma}
\textbf{Proof:}
Since the ring $\bar{R}$ satisfies the Bass stable range condition
Lemma~\ref{congimagest} is a consequence of the familiar argument
that a matrix $g \in \SL_d(K)$ can be reduced by successive applications
of row and column operation to the identity matrix:
Each of these operations is in fact a multiplication by
an elementary matrix from left or right.

\medskip

The following well-known algorithm produces such decomposition (cf. \cite{DV}, Remark 10):

Let $g\in \GL_d(\bar R)$. Let $\bar I$ be the left ideal in $\bar R$
generated by $g_{2d},\dots ,g_{dd}$ in the last column of $g$. The matrix $g$ is left invertible
therefore the left ideal generated by $g_{1d}$ and $\bar I$ coincide with the ring $\bar R$.
The Bass stable range condition implies that there exist $i \in \bar I$ such that
$g_{1d} + i$ is a unit in $\bar R$. This allows us to make last entry
on the first row is an invertible element in $\bar R$ by
multiplying with $d-1$ elementary matrices to the left.
Using an extra multiplication to the left we can make the $d,d$ entry  equal to $1$, next with $d-1$ left and
$d-1$ right multiplications by elementary matrices we can transform the matrix $g$
to an element in $\GL_{d-1}\times 1$ (sitting in the top left corner) .
Thus the reduction form $\GL_d$ to $\GL_{d-1}$ can be done using $3d-2$ elementary matrices,
and by induction using $(3d^2-d-2)/2$ we can transform any matrix to an element in $\GL_1$.

If we begin with an element in $\EL_d$ than it is clear that after applying the above
algorithm the resulting element lies in the intersection of $\GL_1$ with $\EL_d(\bar R)$,
therefore the $1,1$ entry is an invertible element in $\bar R$.
$\square$

\medskip

\noindent
\textbf{Proof of Lemma~\ref{congimage}:}
Let $R = \Mat_t(R')$ where $R'$ is a finite commutative ring.
Start with an element $g \in \EL_d(R)$.
Let $q=\lfloor d/2 \rfloor$ and $p=d-q \geq q$. Using that the ring $\Mat_q(R)$
satisfies the Bass stable range condition we can find GEM
$$
u = \left( \begin{array}{cc}\Id_q & * \\ 0 & \Id_p \end{array}\right)
$$
such that $q\times q$ block in the upper right corner of $ug$ is (left) invertible.
Then there exist a second GEM $u'$ such that lower right $q\times q$ block of $u'ug$ is
equal to $\Id_q$. Using left and right multiplication by two GEMs $u''$ and $v$  we
kill the two off-diagonal blocks of $u''u'ugv$.
Thus we have shown that using $4$ multiplications by generalized elementary matrices we can transform
$g$ to a matrix of type
$$
\left(
\begin{array}{cc}
* & 0 \\ 0 & \Id_q
\end{array}
\right)
$$
where $q=\lfloor d/2 \rfloor$. If we apply the same procedure to the matrix in the upper left corner we can
 transform $g$ to
$$
\left(
\begin{array}{cc}
* & 0 \\ 0 & \Id_{q'}
\end{array}
\right),
$$
where $q' = \lfloor d/2 \rfloor + \lfloor \lceil d/2\rceil/2 \rfloor > d/2$ using $8$ GEMs.

\medskip

It is easy to see that for any $h\in \EL_s(R) = \SL_{st}(R')$ we have the identities
$$
\left(
\begin{array}{cc}
0 & h \\ -h^{-1} & 0
\end{array}
\right)
=
\left(
\begin{array}{cc}
1 & h \\ 0 & 1
\end{array}
\right)
\left(
\begin{array}{cc}
1 & 0 \\ -h^{-1} & 1
\end{array}
\right)
\left(
\begin{array}{cc}
1 & h \\ 0 & 1
\end{array}
\right) \quad \mbox{and}
$$
$$
\left(
\begin{array}{cc}
h & 0 \\ 0 & h^{-1}
\end{array}
\right)
=
\left(
\begin{array}{cc}
1 & 0 \\ 1-h^{-1} & 1
\end{array}
\right)
\left(
\begin{array}{cc}
1 & 1 \\ 0 & 1
\end{array}
\right)
\left(
\begin{array}{cc}
1 & 0 \\ h-1 & 1
\end{array}
\right)
\left(
\begin{array}{cc}
1 & -h^{-1} \\ 0 & 1
\end{array}
\right),
$$
in the group $\EL_{2s}(R) = \SL_{2st}(R')$.

\bigskip

Multiplying the three identities which result from the ones above when we substitute
$h_1$, $-h_2^{-1}$ and $h_1^{-1}h_2^{-1}$ for $h$ in this order gives us that we can write any matrix
$$
\left(
\begin{array}{cc}
[h_1,h_2] & 0 \\ 0 & 1
\end{array}
\right)
=
\left(
\begin{array}{cc}
0& h_1 \\ -h_1^{-1} & 0
\end{array}
\right)
\left(
\begin{array}{cc}
0 & -h_2^{-1} \\ h_2 &0
\end{array}
\right)
\left(
\begin{array}{cc}
h_1^{-1}h_2^{-1} & 0 \\ 0 & h_2h_1
\end{array}
\right)
$$
as a product of $10$ generalized elementary matrices in $\EL_{2s}(R)$.

\medskip

Finally to complete the proof
of Lemma~\ref{congimage} we only need to prove:
\begin{lemma}
\label{commutators}
Let $R'$ be a finite commutative ring. If $s>2$ then any element in $\SL_s(R')$
can be written as a product of $3$ commutators.
\end{lemma}
Without loss of generality we may assume that $R'$ is a finite local commutative ring.
In the case when $R'$ is a field this is a result due to R. C. Thompson
(see \cite{GE} for a survey of the related Ore Conjecture).
\begin{theorem}
\label{ore}
If $K$ is a finite field and $\SL_s(K)$ is perfect (i.e.\ $s>2$ or $|K|>3$)
then every element of it is a commutator.
\end{theorem}

Recall that by a result of Steinberg~\cite{steinberg} the groups $\SL_s(K)$ above
can be generated by 2 elements. We shall need the following basic result:
\begin{lemma}
\label{modules}
Suppose that $\Gamma$ is a group generated by $g_1,..,g_n$, and that $V$ is a right
$\Z \Gamma$-module, such that $V=V(1-\Gamma)$. Then in fact
$$
V= \sum_{i=1}^n V(1-g_i).
$$
\end{lemma}
\textbf{Proof:} Use that
$1-gg_i=g(1-g_i) +(1-g)$ in $\Z \Gamma$ and $g_i,g \in \Gamma$. $\square$
\bigskip

\noindent
\textbf{Proof of Lemma \ref{commutators}:}
Let $I$ be the maximal ideal of $R'$, $K=R' / I$ and define
$H_n:=\ker (\pi_n :\SL_s(R') \to \SL_s(R'/I^n)$
with $H_0=\SL_s(R')$.
It follows that $\bar H:=H_0/H_1\simeq \SL_s(K)$ and
by~\cite{steinberg} it is generated by the images of some $a,b \in H_0$.

On the other hand, each of the quotients
$V_n:=H_n/H_{n+1}$ is an abelian group.
We have that
$$
 [H_1,H_n] \leq H_{n+1}\quad \mathrm{and}\quad [H_0,H_n]=H_n.
$$
The first of these gives that the conjugation action of $H_0$ on $V_n$
makes it into a right $\Z \bar{H}$-module.
The second implies that this module is \emph{perfect}, i.e. that $V_n=V_n(\bar H -1)$ for all $n$.

\medskip

\noindent
\textbf{Claim:} Every $g \in H_1$ can be written as
$$
g=[a,x][b,y]
$$
for appropriate $x,y \in H_1$.

Assuming this for the moment, Theorem \ref{ore} gives that any element $h \in \SL_s(R')$
can be written as $h=[h_1,h_2] g$
with $g \in H_1$ and thus Lemma~\ref{commutators} follows immediately.
\medskip

\noindent
\textbf{Proof of the claim:}
As $H_n=\{1\}$ for large enough $n$ it is enough to prove the following:

\begin{proposition}
\label{series}
Let $g\in H_1$ and $n\in \N$.
There exist elements $x_n, y_n \in H_1$ such that
$$
g= [a,x_n][b,y_n] \quad \textrm{mod } H_n.
$$
\end{proposition}
We shall prove this Proposition by induction starting with the case $n=1$, where it is trivial.
Assuming it is proved for some $n=m-1$ we now prove it for $n=m$:

Without loss of generality, assume that $H_{m+1}=1$, so that $V_{m}=H_m/H_{m+1}$ is central in $H_1$.

Suppose $g= [a,x_n][b,y_n]g'$ for some $g' \in V_m$. We look for
$x_{n+1}=\alpha x_n$, $y_{n+1}=\beta y_n$ with
$\alpha, \beta \in V_m$.

Then by the commutator identities we have
$$
\begin{array}{c@{}c@{}l}
 [a,x_{n+1}][b,y_{n+1}]& {}={}&\displaystyle \rule[10pt]{0pt}{0pt}
  [a,\alpha x_n][b,\beta y_n]=[a,x_n] [a,\alpha]^{x_n} [b,y_n] [b,\beta]^{y_n}= \\
& =& \displaystyle \rule[10pt]{0pt}{0pt}
[a,x_n][b,y_n] [a,\alpha][b,\beta],
\end{array}
$$
the last one because $\alpha,\beta, [a,\alpha],[b,\beta]$ are central in $H_1$.

Therefore we only have to find  $\alpha,\beta \in V_m$ such that
$$
g'= \alpha(1-a)+\beta (1-b)
$$
in the $\bar H$-module $V_m$.
The existence of such $\alpha$ and $\beta$ is a consequence of Lemma~\ref{modules},
recalling that the images of $a$ and $b$ (mod $H_1$)
are generators for $\bar H$.

Proposition~\ref{series} is now proved, and thus Lemma~\ref{commutators} follows.
$\square$

\medskip

Let us finish the proof of Lemma~\ref{congimage}. We have shown that any element $g$
can be transformed to an element in $g' \in \EL_{d'}=\SL_{d't}(R')$ using $8$ GEMs.
If $d't  > 2$ we can apply Lemma~\ref{commutators} to $g'$ and write it
as a product of $3$ commutators, which can be expressed as a product of $10$ GEMs each.
This shows that we can write $g$ as a product of $8+3.10=38$ GEMs.
If $d't \leq 2$ this argument does not work because we ca not
find a copy of $\SL_6(R')$ embedded into $\EL_d(R)$. If $R$ is commutative we can apply
Lemma~\ref{congimagest} and see that every element in $\EL_d(R)$ is a product of
less than $20$ generalized elementary matrices, because after the reduction we end up
with the identity since $\GL_1(R) \cap \EL_d(R) = \{1\}$ .
Finally if $d=3$ and $R=\Mat_2(R')$ we need to modify the above argument slightly
because the group $\SL_2(\F_2)$ is not perfect.
$\square$

\section{Applications}
\label{applications}

\subsection{Unbounded rank expanders}
\label{expanders}

The matrix algebra $\Mat_l(\Z[x_1,\dots,x_{sl}])$ can be generated as a ring
by $s+2$ elements for example the matrices
$$
A = \sum_{i=1}^{l-1} e_{i,i+1} + (-1)^{l-1}e_{l,1} \quad B= e_{1,2} \quad Y_p = \sum_i x_{(p-1)l+i}e_{i,i}, \mbox{ for }p=1,\dots, s.
$$
Therefore we have a surjective homomorphism%
\footnote{There is also a surjective homomorphism $\Z\la x_1,\dots,x_{s+2} \ra \to \Mat_l(\Z[x_1,\dots,x_{sl^2}])$,
but we are unable to use it the obtain estimates of the \KaC\ for $G_{l,sl^2}$ which do not depend on $l$.}
$$
\Z\la x_1,\dots,x_{s+2} \ra \to \Mat_l(\Z[x_1,\dots,x_{sl}]).
$$
This gives a homomorphism
$$
\pi_{l,s}:\EL_3(Z\la x_1,\dots,x_{s+2} \ra) \to
\EL_3(\Mat_l(\Z[x_1,\dots,x_{sl}])) = \SL_{3l}(\Z[x_1,\dots,x_{sl}]).
$$
Let $S_{3,s+2}$ be the standard generating set of $G_{3,s+2}$, consisting of
the $12$ elementary matrices with $\pm 1$ of the diagonal together with
$8(s+2)$ elementary matrices with coefficients the generators $\pm x_i$ next
to the main diagonal. Let $S'_{3l,sl}$ be the generating set of $\SL_{3l}(\Z[x_1,\dots,x_{sl}])$
obtained by applying $\pi_{l,s}$ to $S_{3,s+2}$.

\noindent
\textbf{Proof of Theorem~\ref{tconsforlargerings}:}
Let $\rho :G=\SL_{3l}(\Z[x_1,\dots,x_{sl}]) \to U(\mr{H})$ be a finite unitary
representation and let $v$ be $\epsilon$-almost invariant vector for the set
$S'_{3l,sl}$. Let $H =\ker \rho$ be its kernel. As in the proof of Theorem~\ref{taucom},
let $U$ be the ideal in $\Z[x_1,\dots,x_{sl}]$ generated by the
off diagonal coefficients of the matrices in $H$.

Let $SK_1(R,U;3l)=G(U)/\EL(3l;U)$ where $G(U)$ is the principal congruence subgroup of
$G$ modulo $U$, i.e., the matrices in $G$ congruent to $\Id_{3l}$ modulo $U$.
We have the following diagram
$$
\begin{array}{c@{}c@{}c@{}c@{}c@{}c@{}c@{}c@{}c}
1 & \rightarrow & SK_1(R,U;3l) & \rightarrow & G/\EL(3l;U) & \rightarrow & G(R/U)= \SL_{3l}(\bar{R}) & \rightarrow & 1 \\
  &             &             &             & \downarrow&             &                        &             & \\
  &             &             &             &   G/H     &             &                        &             & \\
  &             &             &             & \downarrow&             &                        &             & \\
  &             &             &             &    1      &             &                        &             & \\
\end{array}
$$
where the row and the column are exact.

We can think of $\SL_{3l}(\bar{R})$ as $\EL_3(\Mat_l(\bar {R}))$ and by Lemma~\ref{congimage} every element
in this group can be written as a product of $38$ generalized elementary matrices in
$\EL_3(\Z\la x_1,\dots,x_{s+2} \ra)$.

As in the proof of Theorem~\ref{taucom} every element in $SK_1(R,U;3l)$ can be written as a product of
$sl+2$ Steinberg symbols which, as in the proof of Theorem~\ref{taucom}, can be written as a product of $13(s+2)$ generalized
elementary matrices.

Therefore every element in $G$ can be written as a product of $64 + 13s$ GEMs and an element
$h$ in the kernel of the representation $\rho$. If we apply Lemma~\ref{abe} to the group
$\EL_3(\Z\la x_1,\dots,x_{s+2}\ra)$ we can see that each of these GEMs moves the
vector $v$ by at most $2M(s+2)\epsilon$.
If $\epsilon$ is small enough than every element in the group $G$ moves $v$ by less than $\sqrt{2}$,
and the representation $\rho$ has an invariant vector. This gives
us a lower bound for the \TauC\ of the form
$$
\TC(\SL_{3l}(\Z[x_1,\dots,x_{sl}]);S'_{3l,sl}) \geq
\frac{1}{\sqrt{2}(64+13s)M(s+2)} \geq \frac{1}{400\left(4+s^{3/2}\right)},
$$
which completes the proof of Theorem~\ref{tconsforlargerings}. $\square$

\medskip

\noindent
\textbf{Proof of Theorem~\ref{UbRExpanders}:}
In the above construction if we put $s=0$ we obtain the \TauC s for the
groups $\SL_{3l}(\Z)$ with respect to the generating sets $S'_{3l}=S'_{3l,0}$, consisting of $28$ elements,
are bounded below by $1/1600$.

If we are only interested in their finite images of the form
$\SL_{3l}(\F_p)$ for $l>3$ then \TauC, which is equal to the \KaC\ because the group is finite,
can be significantly improved.
Let $g \in \SL_{3l}(\F_p) = \SL_3(\Mat_l(\F_p))$.
Using $3$ left multiplications by GEMs we can transform $g$
to a $3\times 3$ block matrix where the last column is trivial,
with an extra $3$ left multiplications by GEMs we can make the
the second column trivial. Finally with one GEM we can transform $g$
to a block diagonal matrix where only the first entry is not trivial.
However this entry is an element in $\SL_l(\F_p)$ and therefore is
a commutator%
\footnote{This is not true if $p=2$ and $l=2$ because $\SL_2(\F_2)$ is not a perfect group. However it is
easy to see that any element in this group is a product of $10$ GEMs in $\EL_3(\Mat_2(\F_2))$.}
 which is expressible as a product of $10$ GEMs, i.e. we have:
$$
\left(\begin{array}{ccc}
* & * & * \\ * & * & * \\ * & * & *
\end{array} \right)
\stackrel{3}{\Longrightarrow}
\left(\begin{array}{ccc}
* & * & 0 \\ * & * & 0 \\ * & * & \Id
\end{array} \right)
\stackrel{3}{\Longrightarrow}
\left(\begin{array}{ccc}
* & 0 & 0 \\ * & \Id & 0 \\ * & 0 & \Id
\end{array} \right)
\stackrel{1}{\Longrightarrow}
$$
$$
\stackrel{1}{\Longrightarrow}
\left(\begin{array}{ccc}
* & 0 & 0 \\ 0 & \Id & 0 \\ 0 & 0 & \Id
\end{array} \right)
\stackrel{10}{\Longrightarrow}
\left(\begin{array}{ccc}
\Id & 0 & 0 \\ 0 & \Id & 0 \\ 0 & 0 & \Id
\end{array} \right)
$$
This shows that every matrix in $\SL_{3l}(\F_p)$ can be written as a product
of $17$ generalized elementary matrices coming form $\SL_3(\Mat_l(\F_p))$.
Therefore if $v$ is $\epsilon$-invariant vector for the set $S'_{3l}$
then $g$ moves $v$ by at most $17\times 2 M(2)\epsilon$, i.e., we have
$$
\TC(\SL_{3l}(\F_p),S'_{3l}) =
\KC(\SL_{3l}(\F_p),S'_{3l}) \geq \frac{\sqrt{2}}{17 \times 2M(2)} \geq \frac{1}{400},
$$
which completes the proof of Theorem~\ref{UbRExpanders}.

\subsection{Luboztky-Wiess Conjectures}
\label{conj}

In~\cite{lubwiess} Luboztky and Wiess conjectured that
for a family of groups the property that their Cayley graphs
are expanders is independant of the generating set.
Recently in~\cite{alonlub}, Alon, Lubotzky and Wigderson
disproved this conjecture using zig-zag products of graphs.
They showed that the Caylay graphs of
$\Gamma_p = \SL_2(\F_p) \ltimes F_2^{p+1}$
are not expanders with respect to some natural generating set, but
that there are expanders with respect to a `random' generating set.
The proof is based on a probabilistic arguments and does not
give explicit generating sets which make the Cayley graphs expanders.

Theorem~\ref{UbRExpanders} together with the observation that
the Cayley graphs of $\SL_{3l}(\F_p)$ are not expanders with respect
to some natural generating sets, provides a natural and explicit
counter example of Conjecture~\ref{conjexp}.

\bigskip

In~\cite{lubwiess}, Luboztky and Wiess also conjectured that an infinite
compact group $K$ can not contain a finitely generated amenable dense subgroup
and a f.g.\ dense subgroup having property \emph{T}. This conjecture
generalizes the observation that a discrete amenable group with
property \emph{T} is finite.
The main motivation for this conjecture is the products of the groups $\SL_n(\F_p)$.
It is known that in
$$
K_n:=\prod_p\SL_n(\F_p)
$$
some natural f.g.\ dense subgroups (for example $\SL_n(\Z)$) have property
\emph{T} and there is no f.g.\ dense amenable subgroup.
On the other hand in
$$
K^p:=\prod_{n}\SL_n(\F_p)
$$
some natural f.g.\ dense subgroups are amenable and there were no known \linebreak[2]
subgroups with property \emph{T} or $\tau$.

\medskip

Using non commutative universal latices, we were
able to construct a \linebreak[4]
residually finite%
\footnote{Note that a discrete residually finite amenable group with
property $\tau$ is also finite.}
f.g.\ subgroup of  $K^p$ which have property $\tau$.
We believe that this group also has property $T$ but we are unable to prove it.
For technical reasons, we are going to work not with subgroups of $K^p$ but with a
similar infinite product.

Let $\widetilde{K^p}$ denote the infinite  product
$$
\widetilde{K^p}:= \prod_{l\geq 4} \SL_{3l}(\F_p).
$$
This is an infinite compact group, because it is an
infinite product of finite, therefore compact, groups.

\medskip

Let $\Gamma_1$ be the subgroup of $\widetilde{K^p}$ generated by
$A$ and $B$, where the components of $A$ in each $\SL_{3l}(\F_p)$
are the matrices $A_{3l} = \sum e_{i,i+1} \pm e_{3l,1}$ and the components of $B$ are
$B_{3l} = \Id_{3l} + e_{1,2}$. Lubotzky and Wiess have shown:
\begin{lemma}[\cite{lubwiess}, Example 4.3]
\label{ame}
The group $\Gamma_1$ is finitely generated dense amenable subgroup
of $\widetilde{K^p}$.
\end{lemma}
\textbf{Proof:}
First we will show that $\Gamma_1$ contains the infinite direct sum
$$
\bigoplus_{l\geq 3} \SL_{3l}(\F_p) \subset \widetilde{K^p},
$$
which implies that $\Gamma_1$ is dense in $\widetilde{K^p}$.

It is enough to show that for each $l\geq 3$, the group $\Gamma_1$ contains the embedding of
$\SL_{3l}(\F_p)$ in $\widetilde{K^p}$ -- the proof uses induction on $l$.
Suppose that we have shown this for all $l<k$.
The direct computation shows that the group commutator of $B$ and $A^{-n}BA^{n}$
has a nontrivial component in $\SL_{3s}(\F_p)$ if and only if $3s \,|\, n-1$ or
$3s\, |\, n+1$. Using the induction hypotheses and the element
$[B,A^{-3k+1}BA^{3k-1}]$ we can see that $\Gamma_1$ contains
element whose only non trivial component lies in the $\SL_{3k}(\F_p)$.
The group $\SL_{3k}(\F_p)$ is quasi-simple and it is normally generated by any of its elements.
This together with the observation that that $\Gamma_1$ projects onto
$\SL_{3k}(\F_p)$ completes the induction step.

\medskip

Let $\Delta_k$ be the subgroup of $\Gamma_1$ generated by $B$, $A^{-1}BA$,\dots,
$A^{-k}BA^k$. It is easy to see that $\Delta_k$ is a finite group, thus
$\Delta= \bigcup_k \Delta_k$ is an increasing union of finite groups and therefore
is amenable. The same argument gives that the group $\widetilde{\Delta} = \bigcup_i A^i\Delta A^{-i}$
is also amenable. By construction $\widetilde{\Delta}$ is the normal closure of $B$ in
$\Gamma_1$. Thus we have the exact sequence
$$
1 \longrightarrow \widetilde{\Delta} \longrightarrow \Gamma_1 \longrightarrow \Z \longrightarrow 1,
$$
i.e., $\Gamma_1$ is extension of an amenable group by an amenable group and it is
also amenable. $\square$

\medskip

There is another way to look at the group $\widetilde{K^p}$
$$
\widetilde{K^p}:= \prod_{l\geq 2} \SL_{3l}(\F_p) =
\prod_{l\geq 2}  \EL_3(\Mat_{l}(\F_p)) =
\EL_3\left(\prod_{l\geq 2}  \Mat_{l}(\F_p)\right).
$$

The ring $R=\prod \Mat_{l}(\F_p)$ contains finitely generated
dense subring $S$, thus group $\Gamma_2=\EL_3(S)$ is a finitely generated
dense subgroup of $\widetilde{K^p}$.
The dense subring $S$ is `almost' isomorphic to the infinite direct sum
$$
\bigoplus \Mat_{l}(\F_p)
$$
and therefore any a finite quotient of $\EL_3(S)$ is `almost' of the form
$\EL_3(S')$ for some finite image $S'$ of the ring $S$. This gives us that
the finite images of $\Gamma_2$ are `almost' finite product of the groups
$\SL_{3l}(\F_p)$ for different $l$.

In the previous section we have shown that the \KaC s, of these groups
with respect to generating sets comming form a generating set of $\EL_3(S)$,
are independent on $l$. This implies
that the group $\EL_3(S)$ has property $\tau$. There are
several technical details, which need to be taken care of, in order to
turn this idea into a proof.

\medskip

Let $S$ be the subring of $R=\prod \Mat_l(\F_p)$ generated by
the elements $A$, $\bar A$, $C$, and $D$, where the components of $A$ in each
matrix algebra $\Mat_l(\F_p)$ are equal to $A_l=\sum e_{i,i+1}\pm e_{l,1}$, the components
of $\bar A$ are $A_l^{-1}$; the components of $C$ and $D$ are
$C_l=e_{1,2}$ and $D_l=e_{2,1}$ respectively.

\begin{lemma}
\label{dense}
The ring $S$ contains the infinite direct sum
$$
\bigoplus_{l\geq 3} \Mat_l(\F_p)
$$
and it is a dense subring of $R$.
\end{lemma}
\textbf{Proof:}
It is enough to show that for any $l$ the ring  $S$ contains the embedding of
the matrix algebra $\Mat_{l}(\F_p)$ in $R$ -- the proof uses induction on $l$.
Suppose that we have shown this for all $l<k$.
The direct computation shows that the ring commutator of $C$ and $\bar A^n C A^{n}$
has nontrivial component in $\Mat_s(\F_p)$ if and only if $s \,|\, n-1$ or
$s\, |\, n+1$. Using the induction hypotheses and the element
$[C,\bar A^{k-1}C A^{k-1}]$ we can see that $S$ contains
an element whose only non trivial component lies in the $\Mat_{k}(\F_p)$.
The algebra $\Mat_{k}(\F_p)$ is simple and the two-sided ideal generated by
any non-zero element coincides with the whole algebra.
This together with the observation that that $S$ projects onto
$\Mat_{k}(\F_p)$ completes the induction step. $\square$

\bigskip

Let $\Gamma_2$ be the subgroup $\EL_3(S) \subset \widetilde{K^p} = \EL_3(R)$.
This group is dense because $S$ is a dense subring of $R$. We also have that
$S$ is a factor ring of $\Z\la x_1,x_2,x_3,x_4\ra$, thus
$\Gamma_2$ is generated by the image of the set $S_{3,4}$, i.e., by
the matrices
$$
e_{i,j}(\pm 1), \quad
e_{i,j}(\pm A), \quad
e_{i,j}(\pm \bar A), \quad
e_{i,j}(\pm C)\quad
\mbox{and} \quad e_{i,j}(\pm D).
$$
This gives us a finitely generated dense subgroup in $\widetilde{K^p}$. Theorem~\ref{amtau} will
follow from Lemmas~\ref{ame} and~\ref{dense} and:
\begin{lemma}
\label{tau}
The group $\Gamma_2$ has property $\tau$.
\end{lemma}
\textbf{Proof:}
Let $\rho$ be a finite representation of $\Gamma_2$ and let $H=\ker \rho$.
As in the proof of Theorem~\ref{taucom} let $U$ be the subset:
$$
U = \{ u \in S \mid e_{i,j}(u) \in H, \mbox{ for some } i\not=j \}.
$$
This is a two sided ideal of finite index in $S$.
By construction $A$ is invertible element in the ring $S$,
therefore there exists $N$ such that
$A^N-1 \in U$. Let $I_N$ be the two-sided ideal in $S$ generated by $A^N-1$.
The inclusion $I_N \subset U$ leads us to the diagram:
$$
\begin{array}{c@{}c@{}c@{}c@{}c@{}c@{}c@{}c@{}c}
1 & \rightarrow & SK_1(S,I_N;3) & \rightarrow & \Gamma_2/\EL_3(I_N) & \rightarrow & \EL_3(S/I_N) & \rightarrow & 1. \\
  &             &             &             & \downarrow&             &                        &             & \\
  &             &             &             & \Gamma_2/H&             &                        &             & \\
\end{array}
$$
In order to finish the proof that $\Gamma_2$ has property $\tau$ it is sufficient to prove that
the groups $\EL_3(S/I_N)$, for all $N$, have uniform bounded elementary generation property, and
that any element in the $K$-group $SK_1(S,I_N;3)$ can be written as a product of small number
of elementary matrices.

Let us consider the factor ring $S/I_N$:

\begin{lemma}
The factor ring $S/I_N$ is isomorphic to
$$
S/I_N \simeq \F_p[C_N] \oplus \bigoplus_{l |N, l\geq 3} \Mat_l(\F_p) ,
$$
where $\F_p[C_N]$ is the group algebra of the cyclic group $C_N$.
\end{lemma}
\textbf{Proof:} 
The inclusion $\bigoplus \Mat_l(\F_p) \subset S$ implies that
$$
S \simeq S^{>2N} \oplus  \bigoplus_{l\geq 3}^{2N} \Mat_l(\F_p)
$$
where $S^{>2N}$ is the the image of $S$ in $\prod_{l>2N}\Mat_l{\F_p}$ .
This isomorphism implies that
$$
\begin{array}{l@{}c@{}l}
S/I_N  & {} \simeq {}&
\displaystyle
S^{>2N}/\left(I_{N}\cap S^{>2N}\right)  \oplus  \bigoplus_{l\geq 3}^{2N} \Mat_l(\F_p)/\left(I_{N}\cap \Mat_l(\F_p)\right) \\
& {}\simeq {} & \displaystyle
S^{>2N}/I_{N}^{>2N} \oplus  \bigoplus_{l\geq 3}^{2N} \Mat_l(\F_p)/I_{N}^{l},
\end{array}
$$
where $I_{N}^{l}$ is the two-sided ideal $\Mat_l(\F_p)$ in generated by $A_l^N-1$
and $I_{N}^{>2N}$ is the two-sided ideal in $S^{>2N}$ generated by the image of $A^N-1$.

The algebra $\Mat_l(\F_p)$ does not have any nontrivial ideals, i.e.,
if $A_l^N-1 \not = 0$ then
$I_{N}^{l} = \Mat_l(\F_p)$.
The last condition is equivalent to $l| N$.
Therefore we have
the isomorphism:
$$
\bigoplus_{l\geq 3}^{2N} \Mat_l(\F_p)/I_{N}^{l}\simeq \bigoplus_{l |N, l\geq 2} \Mat_l(\F_p).
$$

Using the definition of $A$, $\bar A$, $C$ and $D$ it can be
verified that in the ring $S^{>2N}$ we have the identities:
$$
\begin{array}{ccc}
\left[ [C,\bar A C A], \bar A D A \right] = C,
& \quad \quad &
[C,\bar A^{N+1} C A^{N+1}] = 0,
\\
\left[ [\bar A D A,D], \bar A C A \right] = D,
&&
[D,\bar A^{N+1} D A^{N+1}] = 0, \rule[15pt]{0pt}{0pt}
\end{array}
$$
because these identities hold in $\Mat_l(\F_p)$ for all $l>2N$.
These identities imply that
$$
\begin{array}{c@{}c@{}c @{}c@{}c@{}c@{}cc@{}c}
C & {}\equiv{} & \left[ [C,\bar A C A], \bar A D A \right]
  & {}\equiv{} &  \left[ [C,\bar A^{N+1} C A^{N+1}], \bar A D A \right]
  & {}\equiv{} & 0 & \left(\mbox{mod } I_N^{>2N}\right) \\
D & {}\equiv{} & \left[ [\bar A D A,D], \bar A C A \right]
  & {}\equiv{} &  \left[ [\bar A^{N+1} D A^{N+1},D], \bar A C A \right]
  & {}\equiv{} & 0 & \left(\mbox{mod } I_N^{>2N}\right) &, \rule[15pt]{0pt}{0pt}
\end{array}
$$
i.e., $S^{>2N}/I_{N}^{>2N}$ is generated by the images of
$A$ and $\bar A$. Therefore the factor ring is isomorphic
to
$$
\Z\la A,\bar A \ra / (A.\bar A=1, A^N=1, p=0).
$$
The last ring is isomorphic to the
group algebra $\F_p[C_N]$. $\square$

\medskip

Therefore we have
$$
\EL_3(S/I_N) = \EL_3(\F_p[C_N]) \times \prod_{l|N} \EL_{3}(\Mat_l(\F_p)).
$$
By Lemma~\ref{congimage} any element in each term is
a product of $40$ elementary matrices, i.e., we have that
$BE_3(S/I_N) < 40$.

\medskip

It remains to prove the group $SK_1(S,I_N;3)$ is small. Using $K$-theory
we can see that this is an image of $K_2(S/I_N)$. The basic properties of
the functor $K_2$ yield
$$
K_2(S/I_N) = 
K_2\left(\F_p[C_N]\right) \oplus \bigoplus_{l|N} K_2\left(\Mat_l(\F_p)\right).
$$
It is known that $K_2\left(\Mat_l(\F_p)\right) = K_2(\F_p) = 0$, therefore
$$
K_2(S/I_N) = K_2\left(\F_p[C_N]\right).
$$
The ring $\F_p[C_N]$ is a finite commutative ring, generated by
$1$ and another element. Lemma~\ref{kernel} gives that any element in
$K_2\left(\F_p[C_N]\right)$ is a product of $3$ Steinberg symbols%
\footnote{Using the results in~\cite{DS1}, it
can be shown that $K_2\left(\F_p[C_N]\right)$ is a trivial group}
and can be written as a product of $39$ elementary matrices.

\medskip

This shows any element of $\Gamma_2$ can be written as a product of
$80$ elementary matrices modulo an element in the kernel of $\rho$.
Applying Theorem~\ref{t1} gives us that $\Gamma_2$ has property $\tau$
and its \TauC\ is at least
$$
\TC\left(\Gamma_2, \left\{
e_{i,j}(\pm 1), e_{i,j}(\pm A), e_{i,j}(\pm \bar A),  e_{i,j}(\pm C),  e_{i,j}(\pm D)
\big|\, 1\leq\! i\not=j\! \leq 3
\right\}\right) > \frac{1}{2200}.
$$
This completes the proof of Lemma~\ref{tau}. $\square$

As we described above Theorem~\ref{amtau} follows immediately form Lemmas~\ref{ame}, \ref{dense} and~\ref{tau}.

\texttt{\\Martin Kassabov, \\
Cornell University, Ithaca, NY 14853-4201, USA. \\
\emph{e-mail:} kassabov@math.cornell.edu }
\end{document}